\documentclass[aip,
amsmath,amssymb,
preprint,%
numerical,
]{revtex4-1}
\usepackage{graphicx}
\usepackage{dcolumn}
\usepackage{bm}
\usepackage{hyperref}
\usepackage[utf8]{inputenc}
\usepackage[T1]{fontenc}
\usepackage{mathptmx}
\usepackage{etoolbox}
\makeatletter
\def\@email#1#2{%
	\endgroup
	\patchcmd{\titleblock@produce}
	{\frontmatter@RRAPformat}
	{\frontmatter@RRAPformat{\produce@RRAP{*#1\href{mailto:#2}{#2}}}\frontmatter@RRAPformat}
	{}{}
}%
\makeatother
\begin{document}

	\title{Effect of forward scattering and rigid top surface on the onset of phototactic bioconvection in an algal suspension illuminated by both oblique collimated and diffuse irradiation}
	\author{S. K. Rajput}
	\altaffiliation[Email: ]{shubh.iitj@gmail.com}
	\affiliation{ 
		Department of Mathematics, PDPM Indian Institute of Information Technology Design and Manufacturing,
		Jabalpur 482005, India.
	}%
	
	

	\begin{abstract}
	The effect of the rigid top surface on the onset of phototactic bioconvection is investigated using linear stability theory for a suspension of forward-scattering phototactic algae in this article. The suspension is uniformly illuminated by both diffuse and oblique collimated flux. The nature of disturbance of bio-convective instability transits from a stationary (overstable) to an overstable (stationary) state as the forward scattering coefficient varies under fixed parameters. In presence of rigid top surface, the suspension becomes more stable as the forward scattering coefficient increases.
	\end{abstract}
	
	\maketitle

	\section{INTRODUCTION}
		A ubiquitous phenomenon of fluid dynamics which leads to patterns formation in a suspension of living microorganisms (such as algae, bacteria etc.) is defined as Bioconvection~\cite{20platt1961,21pedley1992,22hill2005,23bees2020,24javadi2020}. Platt was the first who introduce the term bioconvection in 1961. In most cases, these living microorganisms have a greater density than the base fluid (water) on a small scale, and they swim collectively in an upward direction. But there are some cases available in a natural environment where there is no need for greater density and to float in an upward direction to form a biological pattern~\cite{21pedley1992}  Non-living microorganism does not show this type of behaviour (pattern formation). Microorganisms show a behavioural response in their swimming direction due to their body structure and external stimuli in the natural environment called taxes.  Gravitaxis is due to gravitational acceleration; gyrotaxis occurs  due to gravitational acceleration and viscous torque where microorganism are bottom heavy; phototaxis is movement of microorganisms towards or away from the illumination source. In this article, we discuss the effect of phototaxis only.    

Studies on bioconvection patterns in algal suspensions have found that the interplay of diffuse and oblique collimated irradiation can affect the flow of fluids and concentration patterns of cells~\cite{1wager1911,3kessler1985,26kessler1997,25kessler1986,4williams2011}. Under intense light condition, the bioconvective patterns may either remain unchanged or be disrupted. The response of algae to moderate or intense light is a critical factor in the changes seen in bioconvection patterns due to light intensity~\cite{6hader1987}. Additionally, the absorption and scattering of light by microorganisms can also play an important role in the changes in bioconvective patterns~\cite{15panda2016}. Algal light scattering can be categorized as isotropic or anisotropic, with the latter being further divided into forward and backward scattering based on cell size, shape, and refractive index. Due to their size, algae mainly scatter light in the forward direction in the visible wavelength range.

The bioconvective system is being studied using the phototaxis model proposed by Panda~\cite{8panda2020}, which involves illuminating a forward-scattering algal suspension with both diffuse and vertical collimated flux. In natural environments, algal suspensions are exposed to oblique collimated flux that can penetrate deeper waters and affect the radiation field, influencing the swimming behavior of algae and light intensity profiles. For realistic models of phototaxis for bioconvective instabilities should include both types of flux (oblique and diffuse). This paper investigates bioconvective instability in presence of both types of flux (oblique and diffuse).

\begin{figure}[!ht]
	\centering
	\includegraphics[width=14cm]{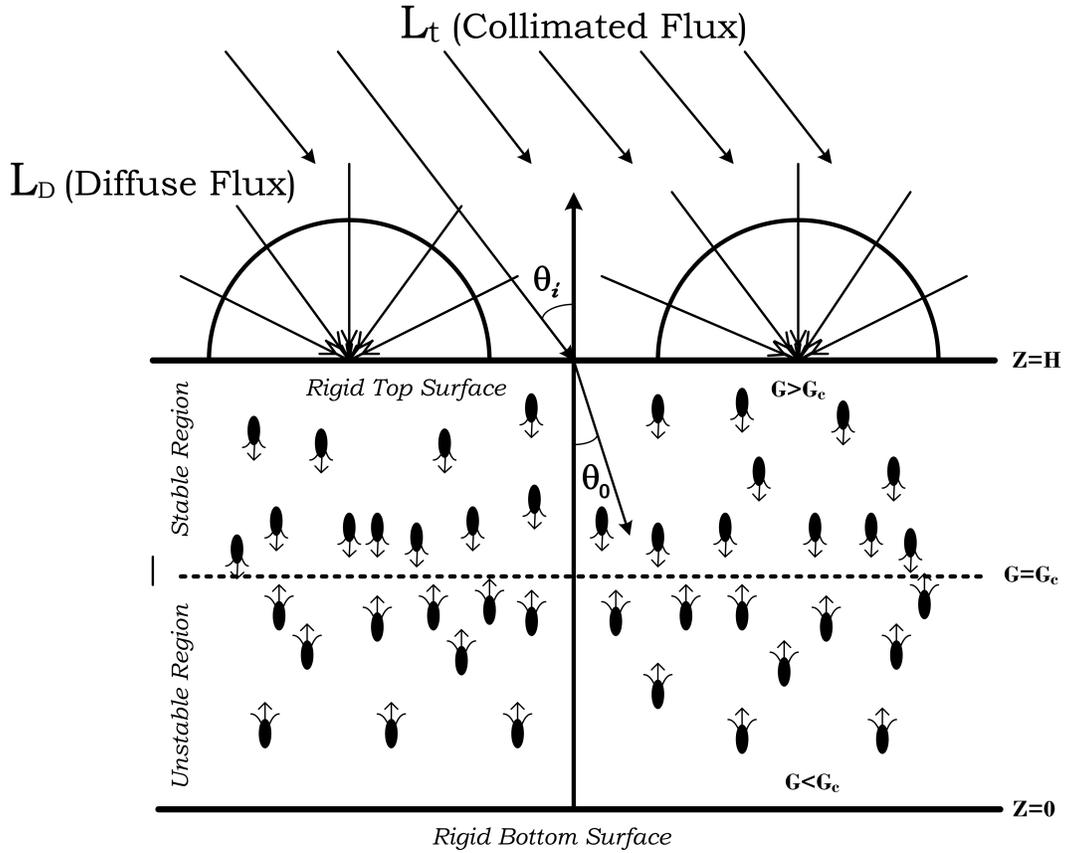}
	\caption{\footnotesize{Accumulation of algae cells as a sublayer at $G=G_c$, where $G_c$ is critical intensity.}}
	\label{fig1}
\end{figure}

This study focuses on the analysis of bioconvective instabilities in a dilute algal suspension. The equilibrium solution is obtained by considering the flow velocity is zero in the suspension. Under this condition, the balance between phototaxis and diffusion leads to a horizontal sublayer where cells accumulate, separating the suspension into two regions. The position of the sublayer is dependent on the critical intensity $G_c$. Above the sublayer $(G > Gc)$, the intensity of light is high enough to suppress cell motion, while below it $(G < Gc)$, cells move upwards towards the sublayer due to phototaxis. This sublayer divides the suspension into two regions based on $G_c$, with the unstable lower region capable of penetrating the stable upper region through fluid motions if the system becomes unstable, a phenomenon referred to as penetrative convection~\cite{9straughan1993}.

	Over the years, researchers have made significant progress in modeling phototactic bioconvection. For instance, Vincent and Hill~\citealp{12vincent1996} studied the onset of phototactic bioconvection and discovered non-stationary and non-oscillatory modes of disturbance where the suspension . Ghorai and Hill~\cite{10ghorai2005} simulated a two-dimensional bioconvective flow pattern using the model of Vincent and Hill~\cite{12vincent1996}, but they neglected scattering by algae. Ghorai et al.~\cite{7ghorai2010} studied the onset of bioconvection in an isotropic scattering suspension and found a steady-state profile with a sublayer located at two different depths due to isotropic scattering by algae for certain parameters. Ghorai and Panda~\cite{13ghorai2013} examined linear stability analysis in a forward scattering algal suspension and reported a transition from a non-oscillatory (non-stationary) to a non-stationary (non-oscillatory) mode at bioconvective instability due to forward scattering. In all these studies, the diffuse flux are not consider as an illuminating source. Panda $et$ $al$.~\cite{11panda2016} investigated the effect of diffuse flux on bioconvcetive instability isotropic scattering algal suspension and found an considerable stabilizing effect due to diffuse flux. After that, Panda~\cite{8panda2020} proposed a model in which he studied the impact of both diffuse and vertical collimated flux and reported that stabilizing effect due to both diffuse collimated flux and forward scattering coefficient. However, these studies did not take into account the effects of oblique collimated flux. Panda $et$ $al.$~\cite{16panda2022} incorporated oblique collimated flux and found that the bioconvective solutions switch from a non-oscillatory state to an overstable state and vice versa at bioconvective instability due to oblique collimated flux. Kumar~\cite{17kumar2022} investigated the impact of oblique collimated flux on bioconvective solutions in an isotropic scattering algal suspension and found that the bioconvection solutions become mostly oscillatory. After that Kumar~\cite{39kumar2023} explored the effect of rigid top and bottom surfaces. In this study, he reprted a considerable stabilizing impact on the bio-convective instabilty caused by rigid surfaces. Panda and Rajput developed a model for a type of algae suspension that scatters light in all directions, and analyzed its stability under illumination by both diffuse and oblique collimated light. They found that the bio-convective solution can transition from non-oscillatory to oscillatory behavior for fixed parameters when subjected to oblique collimated flux. They later extended their analysis to include anisotropic scattering in the suspension and found that the forward scattering coefficient has a stabilizing effect on bio-convective instability. The authors note that further research is needed to examine the impact of both oblique collimated and diffuse light on phototactic bioconvection in a forward scattering algal suspension with rigid non-slip vertical boundaries. The current study aims to investigate the effects of these boundaries on the suspension using a realistic phototaxis model.
	
	The article follows a structured approach. It starts with the mathematical formulation of the problem, followed by obtaining the equilibrium solution and perturbing the base bioconvective governing system by small disturbances. Next, the linear stability problem is derived and solved using numerical methods. Finally, the results of the model are presented and discussed.

	\section{MATHEMATICAL FORMULATION}
	
\begin{figure}[!h]
		\centering
		\includegraphics[width=14cm]{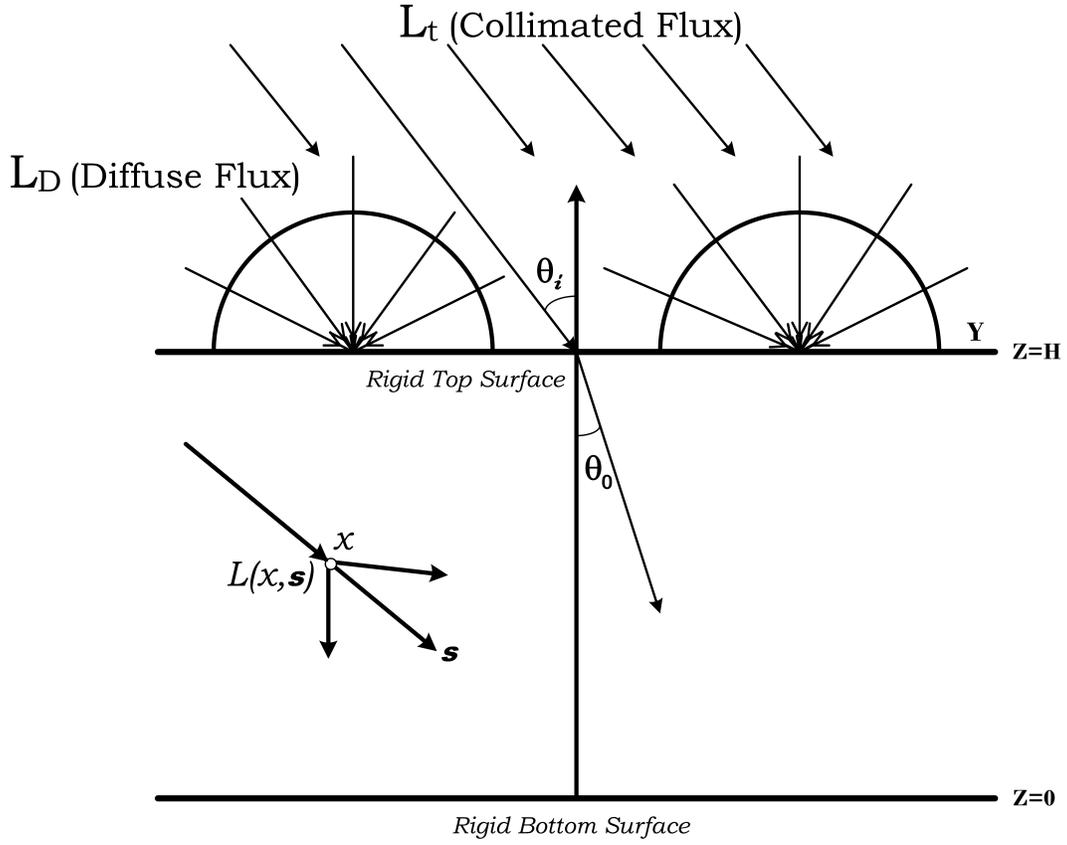}
		\caption{\footnotesize{Geometric configuration of the problem.}}
		\label{fig2}
	\end{figure}

The system being studied is a forward-scattering algal suspension occupying the region between two infinite parallel boundaries in the y-z plane, as depicted in Fig.~\ref{fig2}. The upper surface of the suspension is uniformly illuminated by both oblique collimated and diffuse flux. The algae in the suspension absorb and scatter the incident light in the forward direction, which occurs due to the difference in the refractive index of the algae and water.
	
	\subsection{THE AVERAGE SWIMMING DIRECTION}
	
	The Radiative Transfer Equation (hereafter reffered to as RTE) is used to calculate the light intensity which is given by 
	\begin{equation}\label{1}
		\frac{dL(\boldsymbol{x},\boldsymbol{s})}{ds}+(a+\sigma_s)L(\boldsymbol{x},\boldsymbol{s})=\frac{\sigma_s}{4\pi}\int_{0}^{4\pi}L(\boldsymbol{x},\boldsymbol{s'})\Xi(\boldsymbol{s},\boldsymbol{s'})d\Omega',
	\end{equation}
	where $a,\sigma_s$ and $\Omega'$ are the absorption coefficient, scattering coefficients respectively. $\Xi(\boldsymbol{s},\boldsymbol{s'})$ is the scattering phase function, which provides the angular distribution of the light intensity scattered from $\boldsymbol{s'}$ direction into $\boldsymbol{s}$ direction. Here, we assume that $\Xi(\boldsymbol{s},\boldsymbol{s'})$ is linearly anisotropic with azimuthal symmetry $\Xi(\boldsymbol{s},\boldsymbol{s'})=1+A_1\cos{\theta}\cos{\theta'}$. Here, $A_1$ is the anisotropic coefficient. The anisotropic coefficient defines forward scattering for $0<A_1<1$ and backward scattering for $-1<A_1<0$. If consider the case $A_1=0$, it shows isotropic scattering.

 The light intensity on the top of the suspension is given by
	\begin{equation*}
		L(\boldsymbol{x}_H,\boldsymbol{s})=L_t\delta(\boldsymbol{s}-\boldsymbol{s_0})+\frac{L_D}{\pi}=L_t\delta(\cos{\theta}-\cos{\theta_0})+\frac{L_D}{\pi}, 
	\end{equation*}
	where $\boldsymbol{x}_H=(x,y,H)$ is the location on the top boundary surface. Here, $L_t$ and $L_D$ are the magnitudes of collimated and diffuse irradiation respectively~\cite{8panda2020,15panda2016}.
	  Now consider $a=\alpha n(\boldsymbol{x})$ and $\sigma_s=\beta n(\boldsymbol{ x})$, then the RTE becomes
	\begin{equation}\label{2}
		\frac{dL(\boldsymbol{x},\boldsymbol{s})}{ds}+(\alpha+\beta)nL(\boldsymbol{x},\boldsymbol{s})=\frac{\beta n}{4\pi}\int_{0}^{4\pi}L(\boldsymbol{x},\boldsymbol{s'})(A_1\cos{\theta}\cos{\theta'})d\Omega'.
	\end{equation}
	In the medium, the total intensity at a fixed point $\boldsymbol{x}$ is 
	\begin{equation*}
		G(\boldsymbol{x})=\int_0^{4\pi}L(\boldsymbol{x},\boldsymbol{s})d\Omega,
	\end{equation*}
	and the radiative heat flux is given by
	\begin{equation}\label{3}
		\boldsymbol{q}(\boldsymbol{x})=\int_0^{4\pi}L(\boldsymbol{x},\boldsymbol{s})\boldsymbol{s}d\Omega.
	\end{equation}
	The mean swimming velocity of cell is given by
	\begin{equation*}
		\boldsymbol{W}_c=W_c<\boldsymbol{p}>,
	\end{equation*}
	where $W_c$ is the average cell swimming speed and $<p>$ is cell's mean swimming orientation, which is calculated by
	\begin{equation}\label{4}
		<\boldsymbol{p}>=-M(G)\frac{\boldsymbol{{q}}}{|\boldsymbol{q}|},
	\end{equation}
	where  $M(G), $ is the taxis response function (taxis function), which shows the response of algae cells to light and has the mathematical form such that 
	\begin{equation*}
		M(G)=\left\{\begin{array}{ll}\geq 0, & \mbox{if } G(\boldsymbol{x})\leq G_{c}, \\
			< 0, & \mbox{if }G(\boldsymbol{x})>G_{c}.  \end{array}\right. 
	\end{equation*}
	
	The mean swimming direction becomes zero at the critical light intensity ( $G = G_c)$.
	Generally, the exact functional form of taxis function depends on the species of the microorganisms~\cite{12vincent1996}.

	\subsection{GOVERNING EQUATIONS WITH THE SPECIFIC BOUNDARY CONDITIONS}
	
	Cosider a suspension of phototactic microorganims with average fluid velocity $U$ and $n$ number of algal cells in unit volume. In the algal suspension all cells have constant physical properties like volume $V$ and density $\rho+\Delta\rho$, where $\rho$ is the density of water $(\Delta\rho/\rho<<1)$except for the buoyancy force. Let the suspension be incompressible. Then the sytem of equations of this model is defined as follows:\\
	Continuity equation
	\begin{equation}\label{5}
		\boldsymbol{\nabla}\cdot \boldsymbol{U}=0.
	\end{equation}
	The momentum equation under the Boussinesq approximation 
	\begin{equation}\label{6}
		\rho\frac{D\boldsymbol{U}}{Dt}=-\boldsymbol{\nabla} P+\mu\nabla^2\boldsymbol{U}-ngV\Delta\rho\hat{\boldsymbol{z}},
	\end{equation}
	where  $P$ is dynamic pressure, $\mu$ is the dynamic viscosity of the suspension which is assumed to be that of fluid.\newline 
	Cell conservation equation
	\begin{equation}\label{7}
		\frac{\partial n}{\partial t}=-\boldsymbol{\nabla}\cdot \boldsymbol{B},
	\end{equation}
	where $\boldsymbol{B}$ is the total cell flux which is given by
	\begin{equation}\label{8}
		\boldsymbol{B}=nU+nW_c<\boldsymbol{p}>-\boldsymbol{D}\boldsymbol{\nabla} n.
	\end{equation}
 Here, cell diffusivity $\boldsymbol{D}$ is assumed to be constant, with the result that $\boldsymbol{D} = DI$. Two key assumptions are considered that help to remove Fokker-Plank equation from the giverning equation similar to the Panda.~\cite{8panda2020}
 
Here, lower and upper boundaries are considered rigid and impermeable. Hence, at the boundaries, no fluid flow and no cell movement through the boundaries. The condition of rigid no slip and zero flux are defined as
\begin{equation}\label{9}
	\boldsymbol{U}\times\hat{\boldsymbol{z}}=0\qquad on\quad z=0,H.
\end{equation}
	\begin{equation}\label{10}
		\boldsymbol{B}\cdot\hat{\boldsymbol{z}}=0\qquad on\quad z=0,H.
	\end{equation}
	
We make the assumption that the upper boundary is subjected to oblique and diffuse irradiation uniformly. This leads to a specific set of conditions for the intensity at the boundaries, which are as follows
	\begin{subequations}
		\begin{equation}\label{11a}
			L(x, y, z = 0, \theta, \phi)=L_t\delta(\boldsymbol{s}-\boldsymbol{s_0})+\frac{L_D}{\pi},\quad (\pi/2\leq\theta\leq\pi),
		\end{equation}
		\begin{equation}\label{11b}
			L(x, y, z = 0, \theta, \phi) =0,\quad (0\leq\theta\leq\pi/2).
		\end{equation}
	\end{subequations}

	\subsection{ NON DIMENSIONLIZATION OF THE GOVERNING EQUATIONS}
To express the governing equations without dimensions, certain scales are selected. These include $H$ for length, $H^2/D$ for time, $D/H$ for velocity, $\mu D/H^2$ for pressure, and $\bar{n}$ for concentration. As a result, the governing equations can be represented in a dimensionless form
	\begin{equation}\label{12}
		\boldsymbol{\nabla}\cdot\boldsymbol{U}=0,
	\end{equation}
	\begin{equation}\label{13}
		S_{c}^{-1}\left(\frac{D\boldsymbol{U}}{Dt}\right)=-\nabla P_{e}-Rn\hat{\boldsymbol{z}}+\nabla^{2}\boldsymbol{U},
	\end{equation}
	\begin{equation}\label{14}
		\frac{\partial{n}}{\partial{t}}=-{\boldsymbol{\nabla}}\cdot{\boldsymbol{B}},
	\end{equation}
	where
	\begin{equation}\label{15}
		{\boldsymbol{B}}=n(\boldsymbol{U}+V_{c}<{\boldsymbol{p}}>)-{\boldsymbol{\nabla}}n.
	\end{equation}

In the given equations, the parameter $S_c$ is defined as the ratio of viscosity ($\mu$) and diffusion coefficient $(D)$. The scaled swimming speed is represented by $V_c = W_cH/D$. Additionally, the Rayleigh number $R$ is introduced as a controlling parameter, and is defined as $R=\bar{n}Vg\delta\rho H^3/\mu D$.

The boundary conditions can also be represented in a dimensionless form as
	\begin{equation}\label{16}
	\boldsymbol{U}\times\hat{\boldsymbol{z}}=0\qquad on\quad z=0,1.
\end{equation}
	\begin{equation}\label{17}
		\boldsymbol{B}\cdot\hat{\boldsymbol{z}}=0\qquad on\quad z=0,1.
	\end{equation}	
	
	Nondimensional Radiative Transfer Equation (RTE) is 
	\begin{equation}\label{18}
		\frac{dL(\boldsymbol{x},\boldsymbol{s})}{ds}+\tau_H nL(\boldsymbol{x},\boldsymbol{s})=\frac{\sigma n}{4\pi}\int_{0}^{4\pi}L(\boldsymbol{x},\boldsymbol{s'})(A_1\cos{\theta}\cos{\theta'})d\Omega',
	\end{equation}
	where $\tau_H=(\alpha+\beta)\Bar{n}H$, $\sigma=\beta\Bar{n}H$ are the non-dimensional extinction coefficient and scattering coefficient respectively. The scattering albedo $\omega=\sigma/(\tau_H+\sigma)$ measures the scattering efficiency of microorganisms. In terms of scattering albedo $\omega$, Eq.~(\ref{18}) can be written as
	\begin{equation}\label{19}
		\frac{dL(\boldsymbol{x},\boldsymbol{s})}{ds}+\tau_H nL(\boldsymbol{x},\boldsymbol{s})=\frac{\omega\tau_H n}{4\pi}\int_{0}^{4\pi}L(\boldsymbol{x},\boldsymbol{s'})(A_1\cos{\theta}\cos{\theta'})d\Omega'.
	\end{equation}
The value of the scattering albedo $\omega$ ranges from 0 to 1, where $\omega=1$ represents a medium that purely scatters light, and $\omega=0$ represents a medium that purely absorbs light. The RTE can also be expressed in terms of direction cosine as
	
	\begin{equation}\label{20}
		\xi\frac{dL}{dx}+\eta\frac{dL}{dy}+\nu\frac{dL}{dz}+\tau_H nL(\boldsymbol{x},\boldsymbol{s})=\frac{\omega\tau_H n}{4\pi}\int_{0}^{4\pi}L(\boldsymbol{x},\boldsymbol{s'})(A_1\cos{\theta}\cos{\theta'})d\Omega',
	\end{equation}
In the dimensionless form, the intensity at the boundaries is given by
	\begin{subequations}
		\begin{equation}\label{21a}
			L(x, y, z = 1, \theta, \phi)=L_t\delta(\boldsymbol{s}-\boldsymbol{s_0})+\frac{L_D}{\pi} ,\qquad (\pi/2\leq\theta\leq\pi),
		\end{equation}
		\begin{equation}\label{21b}
			L(x, y, z = 0, \theta, \phi) =0,\qquad (0\leq\theta\leq\pi/2). 
		\end{equation}
	\end{subequations}

	\section{THE BASIC (EQUILIBRIUM) STATE SOLUTION}
	
Equations $(\ref{12})-(\ref{14})$ and $(\ref{20})$ possess a solution at equilibrium that can be expressed as
	
	\begin{equation}\label{22}
		\boldsymbol{U}=0,~~~n=n_s(z)\quad and\quad  L=L_s(z,\theta).
	\end{equation}
	Hence, at the equilibrium state, the total intensity and radiative flux can be expressed as
	
	\begin{equation*}
		G_s=\int_0^{4\pi}L_s(z,\theta)d\Omega,\quad 
		\boldsymbol{q}_s=\int_0^{4\pi}L_s(z,\theta)\boldsymbol{s}d\Omega.
	\end{equation*}
	The $\boldsymbol{q_s}$ has vanishing $x$ and $y$ components because $L_s^d(z,\theta)$ is independent of $\phi$. Therefore, we can express $\boldsymbol{q_s}$ as $\boldsymbol{q}_s=-q_s\hat{\boldsymbol{z}}$, where $q_s=|\boldsymbol{q_s}|$. The equation that governs $L_s$ can be formulated as follows
	\begin{equation}\label{23}
		\frac{dL_s}{dz}+\frac{\tau_H n_sL_s}{\nu}=\frac{\omega\tau_H n_s}{4\pi\nu}(G_s(z)-A_1q_s\nu).
	\end{equation}
	
The equilibrium state intensity can be separated into two components: the collimated part denoted by $L_s^c$ and the diffuse part caused by scattering denoted by $L_s^d$. The equation governing the collimated component $L_s^c$ is  
	
	\begin{equation}\label{24}
		\frac{dL_s^c}{dz}+\frac{\tau_H n_sL_s^c}{\nu}=0,
	\end{equation}
	
	subject to the boundary conditions

	\begin{equation}\label{25}
		L_s^c( 1, \theta) =L_t\delta(\boldsymbol{s}-\boldsymbol{ s}_0),\qquad (\pi/2\leq\theta\leq\pi), 
	\end{equation}
	After calculating the governing equation for $L_s^c$ with boundary condition, we find $L_s^c$
	\begin{equation}\label{26}
		L_s^c=L_t\exp\left(\int_z^1\frac{\tau_H n_s(z')}{\nu}dz'\right)\delta(\boldsymbol{s}-\boldsymbol{s_0}), 
	\end{equation}
	
	and the diffused part is governed by  
	\begin{equation}\label{27}
		\frac{dL_s^d}{dz}+\frac{\tau_H n_sL_s^d}{\nu}=\frac{\omega\tau_H n_s}{4\pi\nu}(G_s(z)-A_1q_s\nu),
	\end{equation}
	subject to the boundary conditions
	\begin{subequations}
		\begin{equation}\label{28a}
			L_s^d( 1, \theta) =\frac{L_D}{\pi},\qquad (\pi/2\leq\theta\leq\pi), 
		\end{equation}
		\begin{equation}\label{28b}
			L_s^d( 0, \theta) =0,\qquad (0\leq\theta\leq\pi/2). 
		\end{equation}
	\end{subequations}

In the basic state, the total intensity $G_s=G_s^c+G_s^d$ is written as
	\begin{equation}\label{29}
		G_s=G_s^c+G_s^d=\int_0^{4\pi}[L_s^c(z,\theta)+L_s^d(z,\theta)d]\Omega=L_t\exp\left(\frac{-\int_z^1\tau_H n_s(z')dz'}{\cos\theta_0}\right)+\int_0^{\pi}L_s^d(z,\theta)d\Omega,
		\end{equation}
	
   Similarly, radiative heat flux in the basic state is defined as
	
	\begin{equation}\label{30}
		\boldsymbol{q_s}=\boldsymbol{q_s}^c+\boldsymbol{q_s}^d=\int_0^{4\pi}\left(L_s^c(z,\theta)+L_s^d(z,\theta)\right)\boldsymbol{s}d\Omega=-L_t(\cos\theta_0)\exp\left(\frac{\int_z^1-\tau_H n_s(z')dz'}{cos\theta_0}\right)\hat{\boldsymbol{z}}+\int_0^{4\pi}L_s^d(z,\theta)\boldsymbol{s}d\Omega.
	\end{equation}
	
	 Now, define a new varible as 
	\begin{equation*}
		\tau=\int_z^1 \tau_H n_s(z')dz',
	\end{equation*}
	
	Eqs.~(\ref{29}) and (\ref{30}) lead to two coupled Fredholm integral equations of the second kind as
	
	\begin{equation}\label{31}
		G_s(\tau) = \frac{\omega}{2}\int_0^{\tau_H} G_s(\tau')E_1(|\tau-\tau'|)d\tau'+e^{-\tau/\cos\theta_0}+2I_DE_2(\tau)+A_1 sgn(\tau-\tau')q_s(\tau')E_2(|\tau-\tau'|),
	\end{equation}

	\begin{equation}\label{32}
	q_s(\tau) = \frac{\omega}{2}\int_0^{\tau_H} A_1q_s(\tau')E_3(|\tau-\tau'|)d\tau'+(\cos\theta_0)e^{-\tau/\cos\theta_0}+2I_DE_3(\tau)+sgn(\tau-\tau')G_s(\tau')E_2(|\tau-\tau'|),
\end{equation}
	
	where $E_n(x)$ is the exponential integral of order $n$ and $sgn(x)$ is the signum function. This coupled FIEs are solved by using method of subtraction of singularity.\par
	
\begin{figure*}[!ht]
	\includegraphics{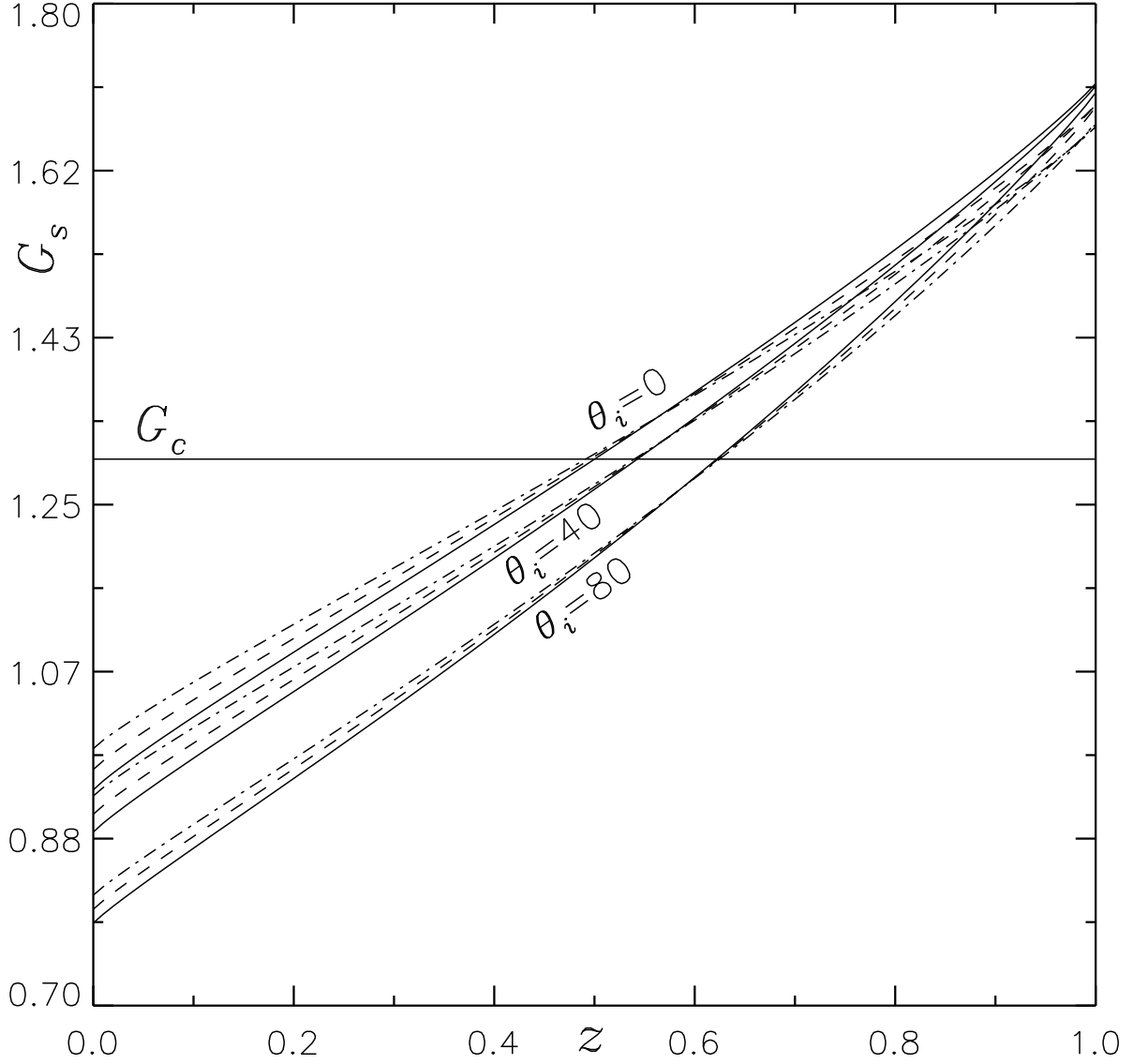}
	\caption{\label{fig3} Variation of total intensity in a uniform suspension by varying forward scattering coefficient $A_1$ from 0 to 0.8 for three different cases $\theta_i=0,40$ and 80. Here, the  governing parameter values $S_c=20,V_c=15,k=0.5,\omega=0.4$ and $L_t=1$ are kept fixed.}
\end{figure*}
	
	The mean swimming direction in the basic state becomes
	
	\begin{equation*}
		<\boldsymbol{p_s}>=-M_s\frac{\boldsymbol{q_s}}{q_s}=M_s\hat{\boldsymbol{z}},
	\end{equation*}
	
	where $M_s=M(G_s).$\par
	In the basic state, the cell concentration $n_s(z)$ satisfy the following equation
	
	\begin{equation}\label{33}
		\frac{dn_s}{dz}-V_cM_sn_s=0,
	\end{equation}
	which is supplemented by the equation
	\begin{equation}\label{34}
		\int_0^1n_s(z)dz=1.
	\end{equation}
	This equation shows the cell conservation relation.
	Eqs.~(\ref{31}) to (\ref{34}) together form a boundary value problem, which can be solved through numerical techniques using the shooting method.

	To illustrate the equilibrium solution, we keep $S_c = 20$, $G_c=1.3$, $J_D = 0.26$, $\tau_H = 0.5$, $\omega = 0.4$, and $(\theta_i)=0,40,80$ fixed. Then, we examine the forward scattering coefficient $A_1$ by varying $A_1$ from 0 to 0.8 on the total intensity and basic equilibrium solution.
	
	Fig.~\ref{fig3} illustrates how the total intensity $G_s$ varies with depth for a uniform suspension as the forward scattering coefficient $A_1$ is increased from 0 to 0.8. As $A_1$ increases, the height of the total intensity $G_s$ at the lower half of the uniform suspension increases, while at the upper half, it decreases. 
	
	\begin{figure*}[!bt]
		\includegraphics{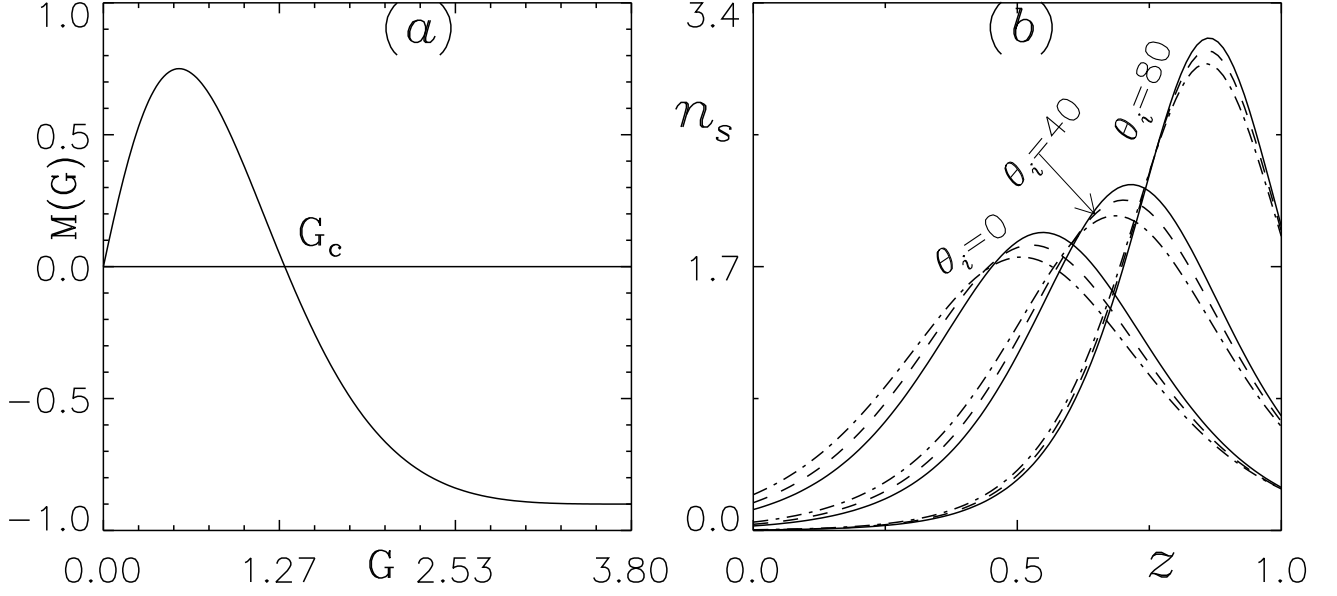}
		\caption{\label{fig4}(a) Taxis response curve for $G_c=1.3$ and, (b) variation in base concentration profile by varying forward scattering coefficient $A_1$ from 0 to 0.8 for three different cases $\theta_i=0,40$ and 80. Here, the  governing parameter values $S_c=20,V_c=15,k=0.5,\omega=0.4$ and $I_t=1$ are kept fixed.}
	\end{figure*}

	Fig.~\ref{4}(a) displays a curve of phototaxis function for light intensity where $G_c = 1.3$ is the critical value. On the other hand, Fig.~\ref{4}(b) presents the influence of the forward scattering coefficient $A_1$ on the sublayer at equilibrium state for the same governing parameters. When $\theta_i = 0$, the sublayer forms at mid-height of the suspension domain under equilibrium state. When $A_1$ is raised from 0 to 0.8, the sublayer location at equilibrium state shifts towards the bottom. Similarly, for $\theta_i = 40$ and 80, the sublayer at equilibrium state is located around three-quarter height and top of the suspension depth, respectively. Increasing A1 from 0 to 0.8 causes the sublayer at equilibrium state to shift towards the bottom in each case.

	\section{Linear stability of the problem}
	To analyze stability, linear perturbation theory is employed, which involves introducing a small perturbation of amplitude $\epsilon<<1$ to the equilibrium state using the following equation

		\begin{align}\label{35}
			\nonumber[\boldsymbol{U},n,L,<\boldsymbol{ p}>]=[0,n_s,L_s,<\boldsymbol{p_s}>]+\epsilon [\boldsymbol{U}_1,n_1,L_1,<\boldsymbol{p}_1>]+\mathcal{O}(\epsilon^2)\\
			=[0,n_s,L_s^c+L_s^d,<\boldsymbol{p}_s>]+\epsilon [\boldsymbol{U}_1,n_1,L_1^c+L_1^d,<\boldsymbol{p}_1>]+\mathcal{O}(\epsilon^2).  
		\end{align}

The perturbed variables are inserted into equations (\ref{12}) to (\ref{14}), and linearization is performed by collecting terms of $o(\epsilon)$ terms about the equilibrium state, which yields
	\begin{equation}\label{36}
		\boldsymbol{\nabla}\cdot \boldsymbol{U}_1=0,
	\end{equation}
	where  $\boldsymbol{U}_1=(U_1,V_1,W_1)$.
	\begin{equation}\label{37}
		S_{c}^{-1}\left(\frac{\partial \boldsymbol{U_1}}{\partial t}\right)+\boldsymbol{\nabla} P_{e}+Rn_1\hat{\boldsymbol{z}}=\nabla^{2}\boldsymbol{ U_1},
	\end{equation}
	\begin{equation}\label{38}
		\frac{\partial{n_1}}{\partial{t}}+V_c\boldsymbol{\nabla}\cdot(<\boldsymbol{p_s}>n_1+<\boldsymbol{p_1}>n_s)+W_1\frac{dn_s}{dz}=\boldsymbol{\nabla}^2n_1.
	\end{equation}
	If $G=G_s+\epsilon G_1+\mathcal{O}(\epsilon^2)=(G_s^c+\epsilon G_1^c)+(G_s^d+\epsilon G_1^d)+\mathcal{O}(\epsilon^2)$, then the steady collimated total intensity is perturbed and after simplification, we get
	\begin{equation}\label{39}
		G_1^c=L_t\exp\left(\frac{-\int_z^1 \tau_H n_s(z')dz'}{\cos\theta_0}\right)\left(\frac{\int_1^z\tau_H n_1 dz'}{\cos\theta_0}\right)
	\end{equation}
	and, $G_1^d$ is given by 
	\begin{equation}\label{40}
		G_1^d=\int_0^{4\pi}L_1^d(\boldsymbol{ x},\boldsymbol{ s})d\Omega.
	\end{equation}
	In the similar manner, we can find perturbed radiative heat flux as 
	\begin{equation}\label{41}
		\boldsymbol{q}_1^c=-I_t(\cos\theta_0)\exp\left(\frac{-\int_z^1 \tau_H n_s(z')dz'}{\cos\theta_0}\right)\left(\frac{\int_1^z\tau_H n_1 dz'}{\cos\theta_0}\right)\hat{z}
	\end{equation}
	and
	\begin{equation}\label{42}
		q_1^d=\int_0^{4\pi}L_1^d(\boldsymbol{ x},\boldsymbol{ s})\boldsymbol{ s}d\Omega.
	\end{equation}
	Now the expression 
	\begin{equation*}
		-M(G_s+\epsilon G_1)\frac{\boldsymbol{q}_s+\epsilon\boldsymbol{q}_1+\mathcal{O}(\epsilon^2)}{|\boldsymbol{q}_s+\epsilon\boldsymbol{q}_1+\mathcal{O}(\epsilon^2)|}-M_s\hat{\boldsymbol{z}},
	\end{equation*}
	gives the perturbed swimming direction on collecting $O(\epsilon)$ terms
	\begin{equation}\label{43}
		<\boldsymbol{p_1}>=G_1\frac{dM_s}{dG}\hat{\boldsymbol{z}}-M_s\frac{\boldsymbol{q_1}^H}{\boldsymbol{q_s}},
	\end{equation}
	where $\boldsymbol{q}_1^H=[\boldsymbol{q}_1^x,\boldsymbol{q}_1^y]$ is the horizontal component of the perturbed radiative flux $\boldsymbol{q}_1$.
	Now substituting the value of $<\boldsymbol{p_1}>$ from  Eq.~$(\ref{46})$ into Eq.~$(\ref{41})$ and simplifying, we get
	\begin{equation}\label{44}
		\frac{\partial{n_1}}{\partial{t}}+V_c\frac{\partial}{\partial z}\left(M_sn_1+n_sG_1\frac{dM_s}{dG}\right)-V_cn_s\frac{M_s}{q_s}\left(\frac{\partial q_1^x}{\partial x}+\frac{\partial q_1^y}{\partial y}\right)+W_1\frac{dn_s}{dz}=\nabla^2n_1.
	\end{equation}
	By taking the double curl and z-component of Eq.$(\ref{37})$, we can eliminate the pressure gradient $P_e$ and the horizontal component of $u_1$. This results in a reduction of Eqs.(\ref{37}), (\ref{38}), and (\ref{44}) to two equations for $W_1$ and $n_1$, which can be decomposed into normal modes as
	\begin{equation}\label{45}
		W_1=\tilde{W}(z)\exp{(\gamma t+i(lx+my))},\quad n_1=\tilde{N}(z)\exp{(\gamma t+i(lx+my))}.  
	\end{equation}
	The governing equation for perturbed intensity $L_1$ can be written as
	\begin{equation}\label{46}
		\xi\frac{\partial L_1}{\partial x}+\eta\frac{\partial L_1}{\partial y}+\nu\frac{\partial L_1}{\partial z}+\tau_H( n_sL_1+n_1L_s)=\frac{\omega\tau_H}{4\pi}(n_sG_1+G_sn_1+A_1\nu(n_sq_1\cdot\hat{z}-q_sn_1)),
	\end{equation}
	subject to the boundary conditions
	\begin{subequations}
		\begin{equation}\label{47a}
			L_1(x, y, z = 1, \xi, \eta, \nu) =0,\qquad (\pi/2\leq\theta\leq\pi,0\leq\phi\leq 2\pi), 
		\end{equation}
		\begin{equation}\label{47b}
			L_1(x, y, z = 0,\xi, \eta, \nu) =0,\qquad (0\leq\theta\leq\pi/2,0\leq\phi\leq 2\pi). 
		\end{equation}
	\end{subequations}
	The Eq. $(\ref{46})$ indicates that $L_1^d$ can be expressed as follows 
	\begin{equation*}
		L_1^d=\Psi^d(z,\xi,\eta,\nu)\exp{(\gamma t+i(lx+my))}. 
	\end{equation*}
	From Eqs.~(\ref{39}) and (\ref{40}), we get
	\begin{equation}\label{48}
		G_1^c=\left[L_t\exp\left(\frac{-\int_z^1 \tau_H n_s(z')dz'}{\cos\theta_0}\right)\left(\frac{\int_1^z\tau_H n_1 dz'}{\cos\theta_0}\right)\right]\exp{(\gamma t+i(lx+my))}=\mathbb{G}^c(z)\exp{(\gamma t+i(lx+my))},
	\end{equation}
	and 
	\begin{equation}\label{49}
		G_1^d=\mathbb{G}^d(z)\exp{(\gamma t+i(lx+my))}= \left(\int_0^{4\pi}\Psi^d(z,\xi,\eta,\nu)d\Omega\right)\exp{(\gamma t+i(lx+my))},
	\end{equation}

	where $\mathbb{G}(z)=\mathbb{G}^c(z)+\mathbb{G}^d(z)$ is the perturbed total intensity. Similarly from Eqs.~(\ref{41}) and (\ref{42}), we have
	\begin{equation*}
		\boldsymbol{q}_1=[q_1^x,q_1^y,q_1^z]=[P(z),Q(z),S(z)]\exp{[\gamma t+i(lx+my)]},
	\end{equation*}
	where 
	\begin{equation*}
		[P(z), Q(z)]=\int_0^{4\pi}[\xi,\eta]\Psi^d(z,\xi,\eta,\nu) d\Omega,
	\end{equation*}

	Note that the $P(z)$ and $Q(z)$ appears due to scattering. On the other hand $S(z)$ has both collimated and diffuse part. The collimated part of $S(z)$ is same as $\boldsymbol{q}_1^c$. Therefore, $S(z)$ is given by

	 \begin{equation*}
		S(z)=\boldsymbol{ q}_1^c+\int_0^{4\pi}\Psi^d(z,\xi,\eta,\nu)\nu d\Omega,
	\end{equation*}
	
	Now, $\Psi^d$ satisfies
	\begin{equation}\label{50}
		\frac{d\Psi^d}{dz}+\frac{(i(l\xi+m\eta)+\tau_H n_s)}{\nu}\Psi^d=\frac{\omega\tau_H}{4\pi\nu}(n_s\mathcal{\mathbb{G}}+G_s\tilde{N}(z)+A_1\nu(n_s S-q_s\tilde{N}(z)))-\frac{\tau_H}{\nu}I_s\tilde{N}(z),
	\end{equation}
	subject to the boundary conditions
	\begin{subequations}
		\begin{equation}\label{51a}
			\Psi^d( 1, \xi, \eta, \nu) =0,\qquad (\pi/2\leq\theta\leq\pi,0\leq\phi\leq 2\pi), 
		\end{equation}
		\begin{equation}\label{51b}
			\Psi^d( 0,\xi, \eta, \nu) =0,\qquad (0\leq\theta\leq\pi/2,0\leq\phi\leq 2\pi). 
		\end{equation}
	\end{subequations}

	The linear stability equations become
	\begin{equation}\label{52}
		\left(\gamma S_c^{-1}+k^2-\frac{d^2}{dz^2}\right)\left( \frac{d^2}{dz^2}-k^2\right)\tilde{W}=Rk^2\tilde{N}(z),
	\end{equation}
	\begin{equation}\label{53}
		\left(\gamma+k^2-\frac{d^2}{dz^2}\right)\tilde{N}(z)+V_c\frac{d}{dz}\left(M_s\tilde{N}(z)+n_s\mathbb{G}\frac{dM_s}{dG}\right)-i\frac{V_cn_sM_s}{q_s}(lP+mQ)=-\frac{dn_s}{dz}\tilde{W}(z),
	\end{equation}
	subject to the boundary conditions
	\begin{equation}\label{54}
		\tilde{W}(z)=\frac{d\tilde{W}(z)}{dz}=\frac{d\tilde{N}(z)}{dz}-V_cM_s\tilde{N}(z)-n_sV_c\mathbb{G}\frac{dM_s}{dG}=0,\quad at\quad z=0,1.
	\end{equation}
	
	Here $k$ represents the non-dimensional wavenumber, which is determined by the square root of the sum of the squares of $l$ and $m$. Eqs.~ (\ref{52})-(\ref{53}) constitute an eigenvalue problem that describes $\gamma$ as a function of various dimensionless parameters, including $V_c$, $\tau_H$, $\omega$, $J_D$, $\theta_i$, $A_1$, $l$, $m$, and $R$. Equation (\ref{53}) can be expressed as 
	\begin{equation}\label{58}
		\Lambda_0(z)+\Lambda_1(z)\int_1^z\tilde{N}(z) dz+(\gamma+k^2+\Lambda_2(z))\tilde{N}(z)+\Lambda_3(z)D\tilde{N}(z)-D^2\tilde{N}(z)=-Dn_s\tilde{W}, 
	\end{equation}
	where
	\begin{subequations}
		\begin{equation}\label{56a}
			\Lambda_0(z)=V_cD\left(n_s\mathbb{G}^d\frac{dM_s}{dG}\right)-\iota\frac{V_cn_sM_s}{q_s}(lP+mQ),
		\end{equation}
		\begin{equation}\label{56b}
			\Lambda_1(z)=\tau_H V_cD\left(n_sG_s^c\frac{dM_s}{dG}\right)
		\end{equation}
		\begin{equation}\label{56c}
			\Lambda_2(z)=2\tau_H V_c n_s G_s^c\frac{dM_s}{dG}+V_c\frac{dM_s}{dM}DG_s^d,
		\end{equation}
		\begin{equation}\label{56d}
			\Lambda_3(z)=V_cM_s.
		\end{equation}
	\end{subequations}
	Introducing the new variable
	\begin{equation}\label{57}
		\tilde{\Theta}(z)=\int_1^z\tilde{N}(z')dz',
	\end{equation}
	the linear stability equations become
	\begin{equation}\label{58}
		\left(\gamma S_c^{-1}+k^2-D^2\right)\left( D^2-k^2\right)\tilde{W}=Rk^2D\tilde{\Theta},
	\end{equation}
	\begin{equation}\label{59}
		\Lambda_0(z)+\Lambda_1(z)\tilde{\Theta}+(\gamma+k^2+\Lambda_2(z))D\tilde{\Theta}+\Lambda_3(z)D^2\tilde{\Theta}-D^3\tilde{\Theta}=-Dn_s\tilde{W}. 
	\end{equation}
	The boundary conditions become,
	
	\begin{equation}\label{60}
		\tilde{W}=D\tilde{W}=D^2\tilde{\Theta}-\Lambda_2(z)D\tilde{\Theta}-\Lambda_3(z)\frac{dM_s}{dG}\mathbb{G}=0,\quad at\quad z=0,1,
	\end{equation}
and the additional boundary condition is
	\begin{equation}\label{61}
		\tilde{\Theta}(z)=0,\quad at\quad z=1.
	\end{equation}

	\section{SOLUTION PROCEDURE}
To solve Eqs. (\ref{58}) and (\ref{59}) and calculate the neutral (marginal) stability curves or the growth rate, $Re(\gamma)$, as a function of R in the (k, R)-plane for a fixed set of other parameters, a fourth-order accurate, finite-difference scheme based on Newton-Raphson-Kantorovich (NRK) iterations~\cite{19cash1980} is utilized. The graph of points for which $Re(\gamma) = 0$ is called a marginal (neutral) stability curve. If the condition $Im(\gamma) = 0$ is satisfied on such a curve, then the bioconvective solution is called stationary (non-oscillatory), and oscillatory solutions exist if $Im(\gamma) \neq 0$. Overstability occurs if the most unstable mode remains on the oscillatory branch of the neutral curve. When an oscillatory solution occurs, a common point $(k_b)$ exists between the stationary and oscillatory branches, and in this instance, the oscillatory branch is the locus of points for which $k \leq k_b$. A particular most unstable mode, i.e. $(k_c, R_c)$, of the neutral curve $R^{(n)}(k)$ (n = 1, 2, 3, ...) is selected, and in this instance, the wavelength of the initial disturbance is calculated as $\lambda_c = 2\pi/k_c$. A bioconvective solution is called mode n if n convection cells can be organized such that one overlies another vertically~\cite{13ghorai2013}. Additionally, the estimated parameters for the proposed problem are the same as previous studies~\cite{7ghorai2010,13ghorai2013,15panda2016,8panda2020,16panda2022}.

	\section{NUMERICAL RESULTS}
	We use a discrete set of fixed parameters to determine the most unstable mode from an initial equilibrium solution. These parameters include $S_c = 20, G_c = 1.3, V_c = 10, 15, 20, \omega = 0.4, J_D = 0.26, 0.48, \theta_i = 0, 40, 80$, and $\tau_H =$ 0.5, 1, while varying $A_1$ as $A_1 =$ 0, 0.4, 0.8. To investigate the influence of forward scattering on bioconvection, we consider discrete values of the angle of incidence ($\theta_i =$ 0, 40, 80). For $\theta_i=0$, the sublayer location at equilibrium state is approximately at mid-height of the domain. As $\theta_i$ increaes to 40 and 80, the sublayer location is shift from mid-height to three-quarter height, and top of the domain, respectively. The width of the unstable region (WUR) refers to the distance from the equilibrium state sublayer location to the bottom of the domain, while the concentration difference in the unstable region (CDUR) is defined as the difference between the maximum concentration and the concentration at the bottom of the suspension. If the WUR or CDUR is higher, it indicates that the suspension is more unstable.
	
	\subsection{WHEN FORWARD SCATTERING IS WEAKER THAN SELF-SHADING  }
	
	The effect of forward scattering on the onset of bioconvection is studied by considering the effectiveness of self-shading compared to scattering, which is achieved by selecting a lower value of the scattering albedo $\omega$. Additionally, the strength of self-shading is varied by selecting different values of the extinction coefficient $\tau_H$, specifically, a high value of $\tau_H=1$ indicates strong self-shading, while a low value of $\tau_H=0.5$ indicates weak self-shading. Throughout the study, a critical intensity of $g_c=1.3$ is used.
	
	\begin{figure*}[!ht]
		\includegraphics{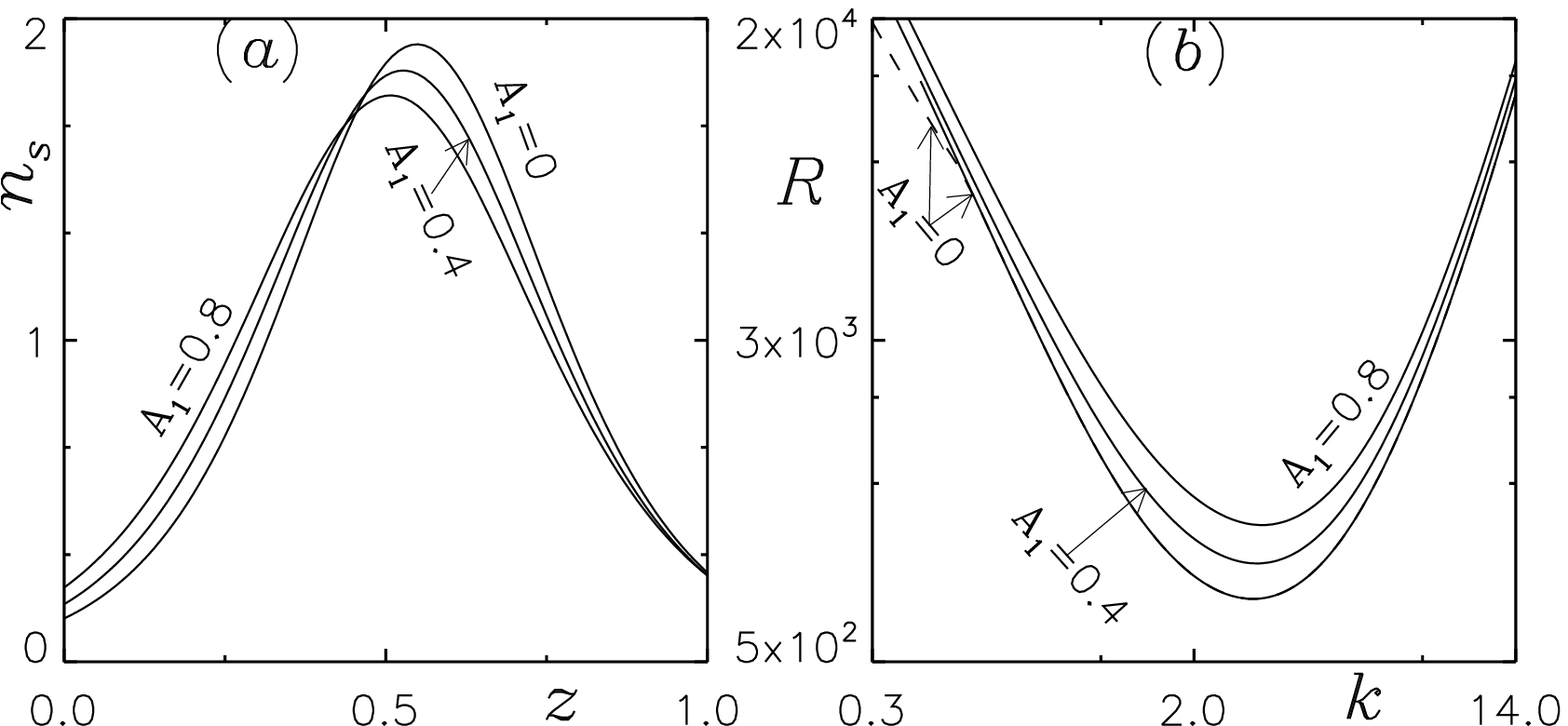}
		\caption{\label{fig5}(a) Basic concentration profile, (b) corresponding marginal stability curves for different values of angle of incidence $\theta_i$. Here, the governing parameter values $S_c=20,V_c=15,k=0.5, I_D=0.26$, and $\omega=0.4$ are kept fixed.}
	\end{figure*}

	\begin{figure*}[!bt]
	    \includegraphics{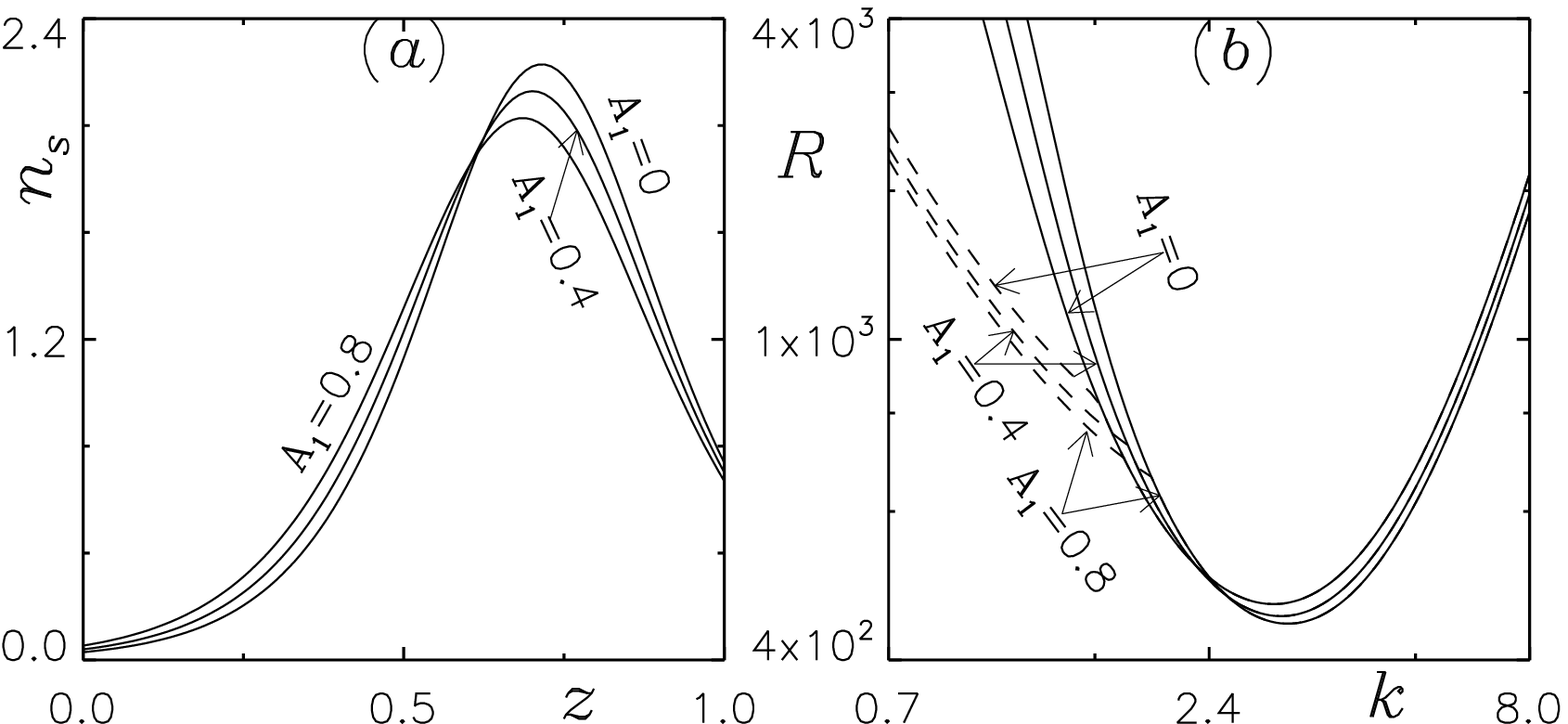}
	    \caption{\label{fig6}(a) Basic concentration profile, (b) corresponding marginal stability curves for different values of angle of incidence $\theta_i$. Here, the governing parameter values $S_c=20,V_c=15,k=0.5, I_D=0.26$, and $\omega=0.4$ are kept fixed.}
    \end{figure*}
	
	\subsubsection{$V_c=15$}
	(\romannumeral 1) When extinction coefficient $\tau_H=0.5$\\
	In this section, the effects of the forward scattering coefficient $A_1$ on bioconvective instability at angle of incidence $\theta_i=0,40,80$ are discussed for a set of fixed parameters $V_c=15,\tau_H=0.5,\omega=0.4$ and $J_D=0.26$. We vary the value of A1 from 0 to 0.8 and study the resulting changes in the sublayer location at equilibrium state, as well as WUR and CDUR.
	
	When $\theta_i=0$, sublayer occurs at the mid-height of the domain for $A_1=0$. As A1 increases, the sublayer at equilibrium state shifts towards the bottom of the domain, resulting in a decrease in both the WUR and CDUR. Consequently, the critical Rayleigh number, which is a measure of the suspension's stability, increases as A1 is increased from 0 to 0.8. As a result, suspension stability increases (see Fig.~\ref{fig5}).
	
	For $\theta_i=40$, the location of the sublayer is at three-quarter height of the domain when $A_1=0$. The location of the sublayer shift towards the midheight of the domain, as $A_1$ increases to 0.4 and 0.8. In this case, WUR and CDUR both decreases as $A_1$ increases. As a result, the critical Rayleigh number increases and suspension becomes more stable similar to the case of $\theta_i=0$ (see Fig.~\ref{fig6}).

		\begin{figure*}[!ht]
		\includegraphics{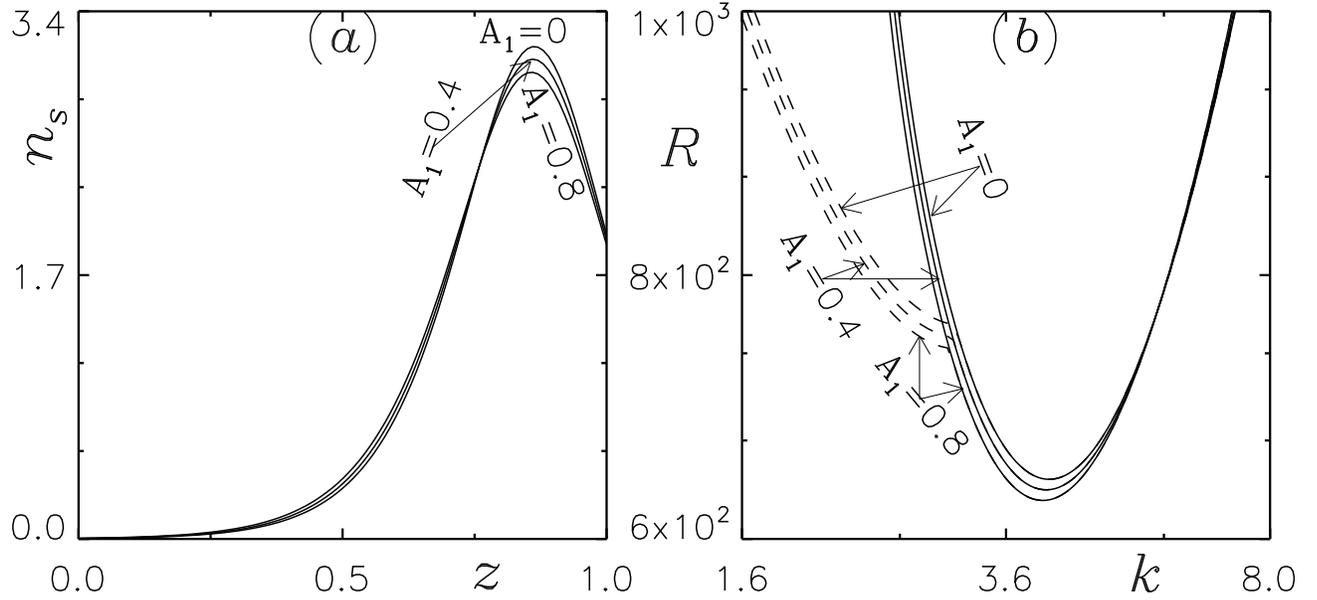}
		\caption{\label{fig7}(a) Basic concentration profile, (b) corresponding marginal stability curves for different values of angle of incidence $\theta_i$. Here, the governing parameter values $S_c=20,V_c=15,k=0.5, I_D=0.26$, and $\omega=0.4$ are kept fixed.}
	\end{figure*}

	For $\theta_i=80$, sublayer formes nearer to the top of the domain. Here, WUR remains same for all cases $A_1=$ 0, 0.4, and 0.8, but CDUR decreases as A1 is increased from 0 to 0.8. In this instance, the decrease in CDUR is more dominant than the change in WUR, resulting in a decrease in critical Rayleigh number and making the suspension less stable (see Fig.~\ref{fig7}).
	
	One more interesting phenomenon (bifurcation of oscillatory branch from the stationary branch) ofthe solution is also observed here. This is observed when $\theta_i=0$ for $A1=0$ and when $\theta_i=40,80$ for all values of $A_1$, but most unstable mode from an initial equilibrium solution is stationary for all values of $\theta_i$.
	
\begin{figure*}[!ht]
	\includegraphics{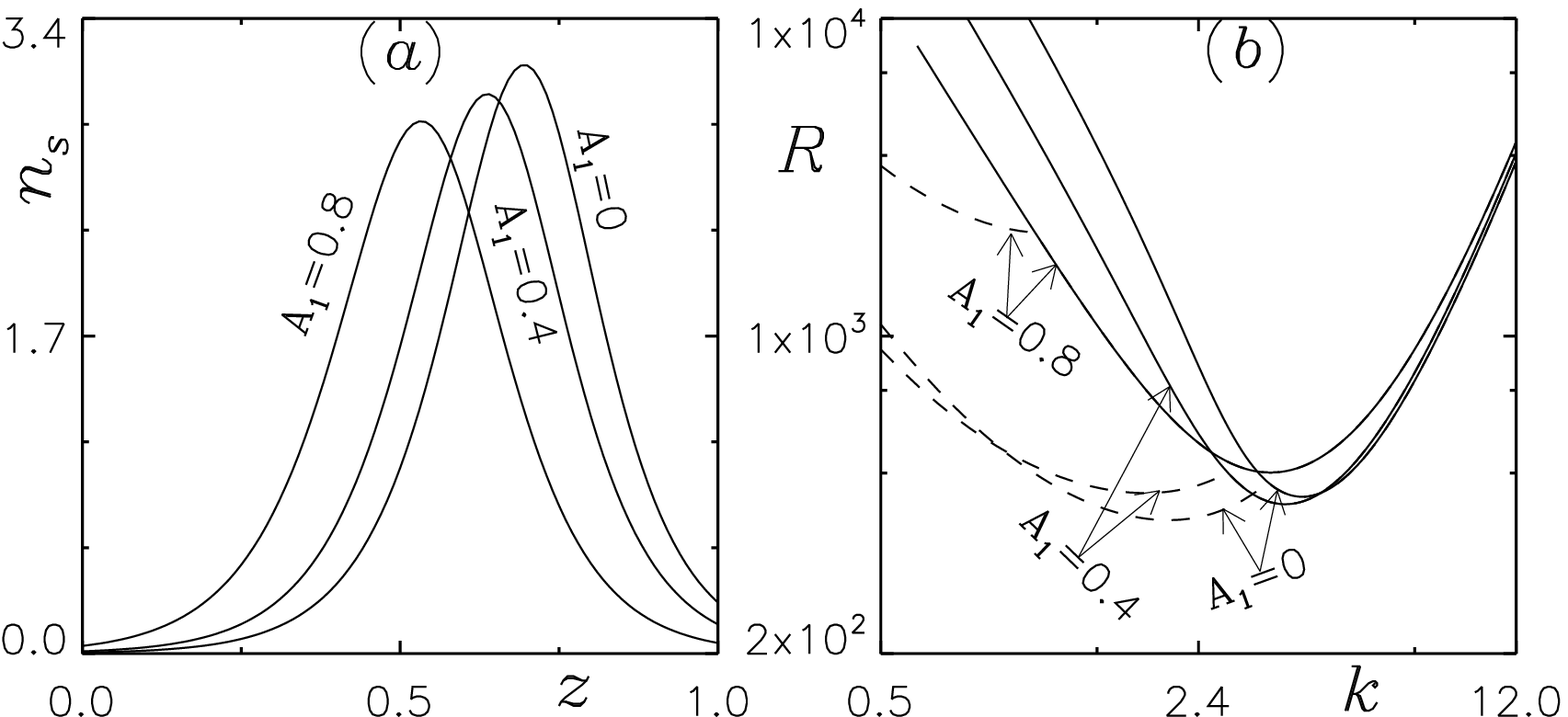}
	\caption{\label{fig8}(a) Basic concentration profile, (b) corresponding marginal stability curves for different values of angle of incidence $\theta_i$. Here, the governing parameter values $S_c=20,V_c=15,k=0.5, I_D=0.26$, and $\omega=0.4$ are kept fixed.}
\end{figure*}

\begin{figure*}[!ht]
	\includegraphics{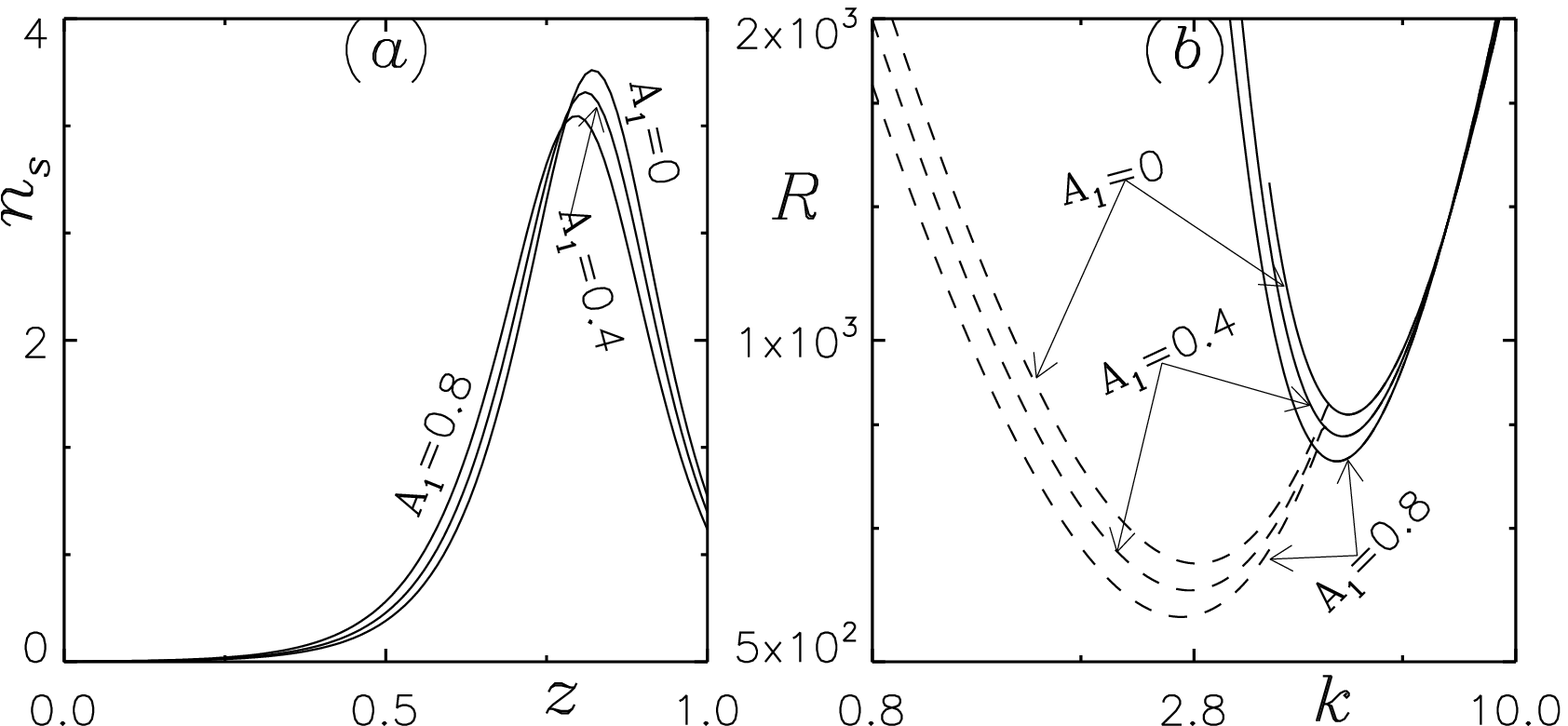}
	\caption{\label{fig9}(a) Basic concentration profile, (b) corresponding marginal stability curves for different values of angle of incidence $\theta_i$. Here, the governing parameter values $S_c=20,V_c=15,k=0.5, I_D=0.26$, and $\omega=0.4$ are kept fixed.}
\end{figure*}

\begin{figure*}[!bt]
	\includegraphics{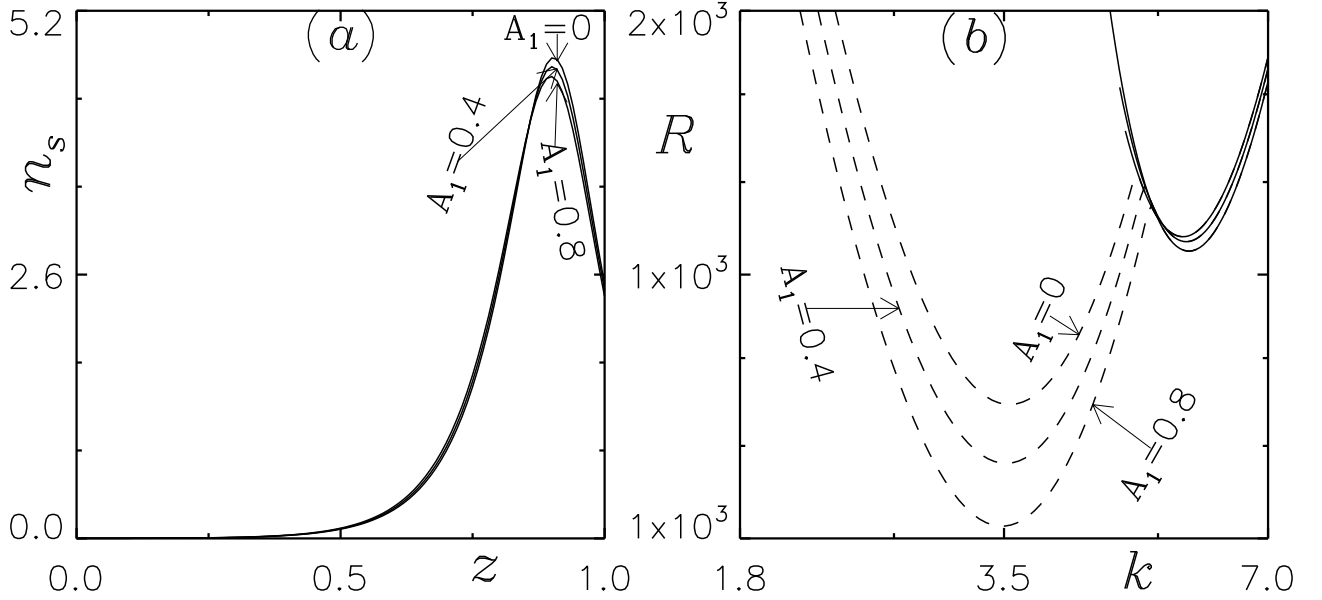}
	\caption{\label{fig10}(a) Basic concentration profile, (b) corresponding marginal stability curves for different values of angle of incidence $\theta_i$. Here, the governing parameter values $S_c=20,V_c=15,k=0.5, I_D=0.26$, and $\omega=0.4$ are kept fixed.}
\end{figure*}

	(\romannumeral 2) When extinction coefficient $\tau_H=1$\\
Let's examine the effect of forward scattering effect where $V_c = 15, \tau_H = 1, \omega = 0.4$, and $J_D = 0.48$ at $\theta_i=0,40,$ and 80.

For the case when $\theta_i=$ 0, the sublayer at equilibrium state forms approximately at a mid-height of the domain, and as $A_1$ increases to 0.4 and 0.8, the location of the sublayer shifts towards the bottom of the domain. Here, WUR and CDUR both decreases as $A_1$ is increased. As a result, critical Rayleigh number is increased and suspension becomes more stable. In this acse, an oscillatory branch bifurcate from the stationary branch of the marginal stabilty curve but only for $A_1=0$, the solution becomes overstable (see Fig.~\ref{fig8}).

Now, for $\theta_i = 40$, When $A_1=0$, the sublayer at equilibrium state forms approximately at around three-quarter height of the domain, and the location of the sublayer at equilibrium state shifts towards the mid-height as $A_1$ increases to 0.4 and 0.8. Here, in WUR a very small decrement is observed with incrent in $A_1$, but CDUR decreases considerably as A1 is increased from 0 to 0.8. In this instance, the decrease in CDUR is more effective than the decrease in WUR, resulting in a decrease in critical Rayleigh number and making the suspension less stable. Here, the most unstable mode in the initial equilibrium solution is overstable for all value of $A_1$.

Finally, for $\theta_i = 80$, the sublayer equilibrium state forms at approximately at the $z=0.9$ of the domain, and the location of the sublayer at equilibrium state shifts towards the three-quarter height of the domain as $A_1$ increases to 0.4 and 0.8. Due to the same region to the case of $\theta_i=40$, the critical Rayleigh number decrease as $A_1$ is increased.  Here, oscillatory branch bifurcates from the stationary branch in all three cases and the most unstable mode in the initial equilibrium solution is overstable. Table~\ref{tab1} summarizes the numerical results for the critical Rayleigh number ($R_c$) and the wavelength ($\lambda_c$) of this section.

	\begin{table}[h]
		\caption{\label{tab1}The quantitative values of bioconvective solutions with increment in $I_D$ for $V_c=15$ are shown in the table, where other parameters are kept fixed.}
		\begin{ruledtabular}
			\begin{tabular}{ccccccccc}
				$V_c$ & $\tau_H$ & $\omega$ & $L_D$ & $\theta_i$ & $A_1$ & $\lambda_c$ & $R_c$ & $Im(\gamma)$ \\
				\hline
				15 & 0.5 & 0.4 & 0.26 & 0 & 0\footnotemark[1] & 2 & 718.98 & 0 \\
				15 & 0.5 & 0.4 & 0.26 & 0 & 0.4 & 2 & 882.36 & 0 \\
				15 & 0.5 & 0.4 & 0.26 & 0 & 0.8 & 1.96 & 1101.18 & 0 \\
				15 & 0.5 & 0.4 & 0.26 & 40 & 0\footnotemark[1] & 1.53 & 452.20 & 0 \\
				15 & 0.5 & 0.4 & 0.26 & 40 & 0.4\footnotemark[1] & 1.57 & 463.92 & 0 \\
				15 & 0.5 & 0.4 & 0.26 & 40 & 0.8\footnotemark[1] & 1.62 & 483.92 & 0 \\
				15 & 0.5 & 0.4 & 0.26 & 80 & 0\footnotemark[1] & 1.2 & 648.90 & 0 \\
				15 & 0.5 & 0.4 & 0.26 & 80 & 0.4\footnotemark[1] & 1.8 & 640.01 & 0 \\
				15 & 0.5 & 0.4 & 0.26 & 80 & 0.8\footnotemark[1] & 1.8 & 631.11 & 0 \\
				
				15 & 1 & 0.4 & 0.48 & 0 & 0 & 2.96\footnotemark[2] & 459.78\footnotemark[2] & 12.58 \\
				15 & 1 & 0.4 & 0.48 & 0 & 0.4\footnotemark[1] & 1.67 & 507.67 & 0 \\
				15 & 1 & 0.4 & 0.48 & 0 & 0.8\footnotemark[1] & 1.79 & 618.54 & 0 \\
				15 & 1 & 0.4 & 0.48 & 40 & 0 & 2.18\footnotemark[2] & 618.01\footnotemark[2] & 20.06 \\
				15 & 1 & 0.4 & 0.48 & 40 & 0.4 & 2.24\footnotemark[2] & 583.28\footnotemark[2] & 19.22 \\
				15 & 1 & 0.4 & 0.48 & 40 & 0.8 & 2.32\footnotemark[2] & 550.74\footnotemark[2] & 18.19 \\
				15 & 1 & 0.4 & 0.48 & 80 & 0 & 1.74\footnotemark[2] & 1190.06\footnotemark[2] & 21.01 \\
				15 & 1 & 0.4 & 0.48 & 80 & 0.4 & 1.76\footnotemark[2] & 1149.62\footnotemark[2] & 21.89 \\
				15 & 1 & 0.4 & 0.48 & 80 & 0.8 & 1.77\footnotemark[2] & 1108.00\footnotemark[2] & 22.53 \\
			\end{tabular}
		\end{ruledtabular}
		\footnotetext[1]{A result indicates that the $R^{(1)}(k)$ branch of the neutral curve is oscillatory.}
		\footnotetext[2]{A result indicates that a smaller solution occurs on the oscillatory branch.}
	\end{table}

	\subsubsection{$V_c=10$ and $V_c=20$}
	The results of the numerical analysis for the bioconvective instability at different values of the forward scattering coefficient for $V_c =$ 10 and 20 show a similar trend to that observed for $V_c =$ 15 in terms of the behavior of the most unstable mode from an equilibrium solution. The details of the numerical results for bioconvective instability are presented in Table~\ref{tab2} and Table~\ref{tab3}.

			\begin{table}[!ht]
			\caption{\label{tab2}The quantitative values of bioconvective solutions with increment in $I_D$ for $V_c=15$ are shown in the table, where other parameters are kept fixed.}
			\begin{ruledtabular}
				\begin{tabular}{ccccccccc}
					$V_c$ & $\tau_H$ & $\omega$ & $L_D$ & $\theta_i$ & $A_1$ & $\lambda_c$ & $R_c$ & $Im(\gamma)$ \\
					\hline
					10 & 0.5 & 0.4 & 0.26 & 0 & 0 & 2.24 & 1335.27 & 0 \\
					10 & 0.5 & 0.4 & 0.26 & 0 & 0.4 & 2.2 & 1577.07 & 0 \\
					10 & 0.5 & 0.4 & 0.26 & 0 & 0.8 & 2.11 & 1884.13 & 0 \\
					10 & 0.5 & 0.4 & 0.26 & 40 & 0 & 2.61 & 619.83 & 0 \\
					10 & 0.5 & 0.4 & 0.26 & 40 & 0.4 & 2.61 & 676.10 & 0 \\
					10 & 0.5 & 0.4 & 0.26 & 40 & 0.8 & 2.61 & 745.77 & 0 \\
					10 & 0.5 & 0.4 & 0.26 & 80 & 0\footnotemark[1] & 2.15 & 457.32 & 0 \\
					10 & 0.5 & 0.4 & 0.26 & 80 & 0.4\footnotemark[1] & 2.2 & 460.12 & 0 \\
					10 & 0.5 & 0.4 & 0.26 & 80 & 0.8 & 2.24 & 463.82 & 0 \\
				
					10 & 1 & 0.4 & 0.48 & 0 & 0\footnotemark[1] & 1.96 & 592.73 & 0 \\
					10 & 1 & 0.4 & 0.48 & 0 & 0.4\footnotemark[1] & 2 & 700.48 & 0 \\
					10 & 1 & 0.4 & 0.48 & 0 & 0.8 & 1.96 & 886.08 & 0 \\
					10 & 1 & 0.4 & 0.48 & 40 & 0\footnotemark[1] & 1.69 & 533.12 & 0 \\
					10 & 1 & 0.4 & 0.48 & 40 & 0.4\footnotemark[1] & 1.72 & 531.64 & 0 \\
					10 & 1 & 0.4 & 0.48 & 40 & 0.8\footnotemark[1] & 1.75 & 534.64 & 0 \\
					10 & 1 & 0.4 & 0.48 & 80 & 0\footnotemark[1] & 1.53 & 695.41 & 0 \\
					10 & 1 & 0.4 & 0.48 & 80 & 0.4\footnotemark[1] & 1.53 & 692.53 & 0 \\
					10 & 1 & 0.4 & 0.48 & 80 & 0.8 & 1.53 & 688.86 & 0 \\
				\end{tabular}
			\end{ruledtabular}
			\footnotetext[1]{A result indicates that the $R^{(1)}(k)$ branch of the neutral curve is oscillatory.}
			\footnotetext[2]{A result indicates that a smaller solution occurs on the oscillatory branch.}
		\end{table}
		
	\begin{table}[h]
		\caption{\label{tab3}The quantitative values of bioconvective solutions with increment in $I_D$ for $V_c=15$ are shown in the table, where other parameters are kept fixed.}
		\begin{ruledtabular}
			\begin{tabular}{ccccccccc}
				$V_c$ & $\tau_H$ & $\omega$ & $L_D$ & $\theta_i$ & $A_1$ & $\lambda_c$ & $R_c$ & $Im(\gamma)$ \\
				\hline
				20 & 0.5 & 0.4 & 0.26 & 0 & 0\footnotemark[1] & 2 & 478.48 & 0 \\
				20 & 0.5 & 0.4 & 0.26 & 0 & 0.4\footnotemark[1] & 2 & 577.36 & 0 \\
				20 & 0.5 & 0.4 & 0.26 & 0 & 0.8\footnotemark[1] & 1.96 & 753.15 & 0 \\
				20 & 0.5 & 0.4 & 0.26 & 40 & 0\footnotemark[1] & 1.53 & 552.23 & 0 \\
				20 & 0.5 & 0.4 & 0.26 & 40 & 0.4\footnotemark[1] & 1.57 & 536.7 & 0 \\
				20 & 0.5 & 0.4 & 0.26 & 40 & 0.8\footnotemark[1] & 1.62 & 524.43 & 0 \\
				20 & 0.5 & 0.4 & 0.26 & 80 & 0\footnotemark[1] & 1.2 & 1049.65 & 0 \\
				20 & 0.5 & 0.4 & 0.26 & 80 & 0.4 & 1.8\footnotemark[2] & 1017\footnotemark[2] & 16.06 \\
				20 & 0.5 & 0.4 & 0.26 & 80 & 0.8 & 1.8\footnotemark[2] & 984.16\footnotemark[2] & 15.63 \\
				
				20 & 1 & 0.4 & 0.48 & 0 & 0 & 2.67\footnotemark[2] & 423.32\footnotemark[2] & 21.58 \\
				20 & 1 & 0.4 & 0.48 & 0 & 0.4 & 3.04\footnotemark[2] & 379.34\footnotemark[2] & 17.56 \\
				20 & 1 & 0.4 & 0.48 & 0 & 0.8 & 2.72\footnotemark[2] & 718.49\footnotemark[2] & 12.04 \\
				20 & 1 & 0.4 & 0.48 & 40 & 0 & 1.93\footnotemark[2] & 820.65\footnotemark[2] & 35.14 \\
				20 & 1 & 0.4 & 0.48 & 40 & 0.4 & 1.96\footnotemark[2] & 750.14\footnotemark[2] & 33.71 \\
				20 & 1 & 0.4 & 0.48 & 40 & 0.8 & 2.05\footnotemark[2] & 680.13\footnotemark[2] & 31.83 \\
				20 & 1 & 0.4 & 0.48 & 80 & 0 & 1.49\footnotemark[2] & 1854.42\footnotemark[2] & 38.81 \\
				20 & 1 & 0.4 & 0.48 & 80 & 0.4 & 1.52\footnotemark[2] & 1779.89\footnotemark[2] & 40.09 \\
				20 & 1 & 0.4 & 0.48 & 80 & 0.8 & 1.52\footnotemark[2] & 1702.46\footnotemark[2] & 40.98 \\
			\end{tabular}
		\end{ruledtabular}
		\footnotetext[1]{A result indicates that the $R^{(1)}(k)$ branch of the neutral curve is oscillatory.}
		\footnotetext[2]{A result indicates that a smaller solution occurs on the oscillatory branch.}
	\end{table}

	\section{Conclusion}
The proposed phototaxis model investigates the impact of rigid top surface on the onset of light induced bioconvection in an anisotropic (forward) scattering algal suspension illuminated by both diffuse and obliqe (not vertical) collimated flux. A linear anisotropic scattering coefficent is used in this analysis. An initial equilibrium solution at the bioconvective instability is also investigated by using the linear analysis of the same suspenaion.

%

 When the forward scattering coefficient $A_1$ is increased in a uniform suspension, the total intensity decreases in the upper half and increases in the lower half of the suspension at the equilibrium state. Furthermore, the critical value of total intensity decrease (increase) in the lower (upper) half of the suspension as $A_1$ is increased. As $A_1$ is increased, the sublayer position shift towards the bottom of the domain. On the other hand, concentration in the sublayer also increases as the $A_1$ is increased.
	
When the sublayer at equilibrium state forms near the mid-height of the domain, the critical wavelength decreases and critical Rayleigh number increases with an increase in the forward scattering coefficient $A_1$. When sublayer forms nearer to the top of the domain, as $A_1$ is increased, the critical Rayleigh number and critical wavelength both decreases. On the other hand, when sublayer occurs around at three-quarter height of the domain, the critical rayleigh number and critical wavelength both increases for weak self-shading suspension but for strong self-shading suspension, the critical wavelength increases and critical Rayleigh number decreases as $A_1$ is increased. 

It seems that the behavior of the bioconvective solutions depends on the self-shading effect in the suspension. For weak self-shading suspension, an oscillatory branch bifurcate from the stationary branch of the marginal stability curve but most unstable solution occurs on the stationary branch. Therefore, almost every position of the sulayer solution remains stationary. But, for strong self-shading suspension, the bio-convective solution convert into overstable solution for almost all locations for the sublayer at equilibrium state.
	
	
	\begin{acknowledgments}
		The author gratefully acknowledges the Ministry of Education (Government of India) for the ﬁnancial support via GATE fellowship (Registration No. MA19S43047204). 
	\end{acknowledgments}
	
	\section*{Data Availability}
	The data that support the plots within this paper and
	other findings of this study are  available
	within the article.
	\nocite{*}
	\bibliography{main}

\providecommand{\noopsort}[1]{}\providecommand{\singleletter}[1]{#1}%
\begin{thebibliography}{39}%
\makeatletter
\providecommand \@ifxundefined [1]{%
 \@ifx{#1\undefined}
}%
\providecommand \@ifnum [1]{%
 \ifnum #1\expandafter \@firstoftwo
 \else \expandafter \@secondoftwo
 \fi
}%
\providecommand \@ifx [1]{%
 \ifx #1\expandafter \@firstoftwo
 \else \expandafter \@secondoftwo
 \fi
}%
\providecommand \natexlab [1]{#1}%
\providecommand \enquote  [1]{``#1''}%
\providecommand \bibnamefont  [1]{#1}%
\providecommand \bibfnamefont [1]{#1}%
\providecommand \citenamefont [1]{#1}%
\providecommand \href@noop [0]{\@secondoftwo}%
\providecommand \href [0]{\begingroup \@sanitize@url \@href}%
\providecommand \@href[1]{\@@startlink{#1}\@@href}%
\providecommand \@@href[1]{\endgroup#1\@@endlink}%
\providecommand \@sanitize@url [0]{\catcode `\\12\catcode `\$12\catcode
  `\&12\catcode `\#12\catcode `\^12\catcode `\_12\catcode `\%12\relax}%
\providecommand \@@startlink[1]{}%
\providecommand \@@endlink[0]{}%
\providecommand \url  [0]{\begingroup\@sanitize@url \@url }%
\providecommand \@url [1]{\endgroup\@href {#1}{\urlprefix }}%
\providecommand \urlprefix  [0]{URL }%
\providecommand \Eprint [0]{\href }%
\providecommand \doibase [0]{http://dx.doi.org/}%
\providecommand \selectlanguage [0]{\@gobble}%
\providecommand \bibinfo  [0]{\@secondoftwo}%
\providecommand \bibfield  [0]{\@secondoftwo}%
\providecommand \translation [1]{[#1]}%
\providecommand \BibitemOpen [0]{}%
\providecommand \bibitemStop [0]{}%
\providecommand \bibitemNoStop [0]{.\EOS\space}%
\providecommand \EOS [0]{\spacefactor3000\relax}%
\providecommand \BibitemShut  [1]{\csname bibitem#1\endcsname}%
\let\auto@bib@innerbib\@empty
\bibitem [{\citenamefont {Platt}(1961)}]{20platt1961}%
  \BibitemOpen
  \bibfield  {author} {\bibinfo {author} {\bibfnamefont {J.~R.}\ \bibnamefont
  {Platt}},\ }\bibfield  {title} {\enquote {\bibinfo {title} {" bioconvection
  patterns" in cultures of free-swimming organisms},}\ }\href@noop {}
  {\bibfield  {journal} {\bibinfo  {journal} {Science}\ }\textbf {\bibinfo
  {volume} {133}},\ \bibinfo {pages} {1766--1767} (\bibinfo {year}
  {1961})}\BibitemShut {NoStop}%
\bibitem [{\citenamefont {Pedley}\ and\ \citenamefont
  {Kessler}(1992)}]{21pedley1992}%
  \BibitemOpen
  \bibfield  {author} {\bibinfo {author} {\bibfnamefont {T.~J.}\ \bibnamefont
  {Pedley}}\ and\ \bibinfo {author} {\bibfnamefont {J.~O.}\ \bibnamefont
  {Kessler}},\ }\bibfield  {title} {\enquote {\bibinfo {title} {Hydrodynamic
  phenomena in suspensions of swimming microorganisms},}\ }\href@noop {}
  {\bibfield  {journal} {\bibinfo  {journal} {Annual Review of Fluid
  Mechanics}\ }\textbf {\bibinfo {volume} {24}},\ \bibinfo {pages} {313--358}
  (\bibinfo {year} {1992})}\BibitemShut {NoStop}%
\bibitem [{\citenamefont {Hill}\ and\ \citenamefont
  {Pedley}(2005)}]{22hill2005}%
  \BibitemOpen
  \bibfield  {author} {\bibinfo {author} {\bibfnamefont {N.~A.}\ \bibnamefont
  {Hill}}\ and\ \bibinfo {author} {\bibfnamefont {T.~J.}\ \bibnamefont
  {Pedley}},\ }\bibfield  {title} {\enquote {\bibinfo {title}
  {Bioconvection},}\ }\href@noop {} {\bibfield  {journal} {\bibinfo  {journal}
  {Fluid Dynamics Research}\ }\textbf {\bibinfo {volume} {37}},\ \bibinfo
  {pages} {1} (\bibinfo {year} {2005})}\BibitemShut {NoStop}%
\bibitem [{\citenamefont {Bees}(2020)}]{23bees2020}%
  \BibitemOpen
  \bibfield  {author} {\bibinfo {author} {\bibfnamefont {M.~A.}\ \bibnamefont
  {Bees}},\ }\bibfield  {title} {\enquote {\bibinfo {title} {Advances in
  bioconvection},}\ }\href@noop {} {\bibfield  {journal} {\bibinfo  {journal}
  {Annual Review of Fluid Mechanics}\ }\textbf {\bibinfo {volume} {52}},\
  \bibinfo {pages} {449--476} (\bibinfo {year} {2020})}\BibitemShut {NoStop}%
\bibitem [{\citenamefont {Javadi}\ \emph {et~al.}(2020)\citenamefont {Javadi},
  \citenamefont {Arrieta}, \citenamefont {Tuval},\ and\ \citenamefont
  {Polin}}]{24javadi2020}%
  \BibitemOpen
  \bibfield  {author} {\bibinfo {author} {\bibfnamefont {A.}~\bibnamefont
  {Javadi}}, \bibinfo {author} {\bibfnamefont {J.}~\bibnamefont {Arrieta}},
  \bibinfo {author} {\bibfnamefont {I.}~\bibnamefont {Tuval}}, \ and\ \bibinfo
  {author} {\bibfnamefont {M.}~\bibnamefont {Polin}},\ }\bibfield  {title}
  {\enquote {\bibinfo {title} {Photo-bioconvection: towards light control of
  flows in active suspensions},}\ }\href@noop {} {\bibfield  {journal}
  {\bibinfo  {journal} {Philosophical Transactions of the Royal Society A}\
  }\textbf {\bibinfo {volume} {378}},\ \bibinfo {pages} {20190523} (\bibinfo
  {year} {2020})}\BibitemShut {NoStop}%
\bibitem [{\citenamefont {Wager}(1911)}]{1wager1911}%
  \BibitemOpen
  \bibfield  {author} {\bibinfo {author} {\bibfnamefont {H.}~\bibnamefont
  {Wager}},\ }\bibfield  {title} {\enquote {\bibinfo {title} {Vii. on the
  effect of gravity upon the movements and aggregation of euglena viridis,
  ehrb., and other micro-organisms},}\ }\href@noop {} {\bibfield  {journal}
  {\bibinfo  {journal} {Philosophical Transactions of the Royal Society of
  London. Series B, Containing Papers of a Biological Character}\ }\textbf
  {\bibinfo {volume} {201}},\ \bibinfo {pages} {333--390} (\bibinfo {year}
  {1911})}\BibitemShut {NoStop}%
\bibitem [{\citenamefont {Kessler}(1985)}]{3kessler1985}%
  \BibitemOpen
  \bibfield  {author} {\bibinfo {author} {\bibfnamefont {J.~O.}\ \bibnamefont
  {Kessler}},\ }\bibfield  {title} {\enquote {\bibinfo {title} {Co-operative
  and concentrative phenomena of swimming micro-organisms},}\ }\href@noop {}
  {\bibfield  {journal} {\bibinfo  {journal} {Contemporary Physics}\ }\textbf
  {\bibinfo {volume} {26}},\ \bibinfo {pages} {147--166} (\bibinfo {year}
  {1985})}\BibitemShut {NoStop}%
\bibitem [{\citenamefont {Kessler}\ and\ \citenamefont
  {Hill}(1997)}]{26kessler1997}%
  \BibitemOpen
  \bibfield  {author} {\bibinfo {author} {\bibfnamefont {J.~O.}\ \bibnamefont
  {Kessler}}\ and\ \bibinfo {author} {\bibfnamefont {N.~A.}\ \bibnamefont
  {Hill}},\ }\bibfield  {title} {\enquote {\bibinfo {title} {Complementarity of
  physics, biology and geometry in the dynamics of swimming micro-organisms},}\
  }in\ \href@noop {} {\emph {\bibinfo {booktitle} {Physics of biological
  systems}}}\ (\bibinfo  {publisher} {Springer},\ \bibinfo {year} {1997})\ pp.\
  \bibinfo {pages} {325--340}\BibitemShut {NoStop}%
\bibitem [{\citenamefont {Kessler}(1986)}]{25kessler1986}%
  \BibitemOpen
  \bibfield  {author} {\bibinfo {author} {\bibfnamefont {J.~O.}\ \bibnamefont
  {Kessler}},\ }\bibfield  {title} {\enquote {\bibinfo {title} {The external
  dynamics of swimming micro-organisms},}\ }\href@noop {} {\bibfield  {journal}
  {\bibinfo  {journal} {Progress in phycological research}\ }\textbf {\bibinfo
  {volume} {4}},\ \bibinfo {pages} {258--307} (\bibinfo {year}
  {1986})}\BibitemShut {NoStop}%
\bibitem [{\citenamefont {Williams}\ and\ \citenamefont
  {Bees}(2011)}]{4williams2011}%
  \BibitemOpen
  \bibfield  {author} {\bibinfo {author} {\bibfnamefont {C.~R.}\ \bibnamefont
  {Williams}}\ and\ \bibinfo {author} {\bibfnamefont {M.~A.}\ \bibnamefont
  {Bees}},\ }\bibfield  {title} {\enquote {\bibinfo {title} {A tale of three
  taxes: photo-gyro-gravitactic bioconvection},}\ }\href@noop {} {\bibfield
  {journal} {\bibinfo  {journal} {Journal of Experimental Biology}\ }\textbf
  {\bibinfo {volume} {214}},\ \bibinfo {pages} {2398--2408} (\bibinfo {year}
  {2011})}\BibitemShut {NoStop}%
\bibitem [{\citenamefont {H{\"a}der}(1987)}]{6hader1987}%
  \BibitemOpen
  \bibfield  {author} {\bibinfo {author} {\bibfnamefont {D.-P.}\ \bibnamefont
  {H{\"a}der}},\ }\bibfield  {title} {\enquote {\bibinfo {title} {Polarotaxis,
  gravitaxis and vertical phototaxis in the green flagellate, euglena
  gracilis},}\ }\href@noop {} {\bibfield  {journal} {\bibinfo  {journal}
  {Archives of microbiology}\ }\textbf {\bibinfo {volume} {147}},\ \bibinfo
  {pages} {179--183} (\bibinfo {year} {1987})}\BibitemShut {NoStop}%
\bibitem [{\citenamefont {Panda}\ \emph {et~al.}(2016)\citenamefont {Panda},
  \citenamefont {Singh}, \citenamefont {Mishra},\ and\ \citenamefont
  {Mohanty}}]{15panda2016}%
  \BibitemOpen
  \bibfield  {author} {\bibinfo {author} {\bibfnamefont {M.~K.}\ \bibnamefont
  {Panda}}, \bibinfo {author} {\bibfnamefont {R.}~\bibnamefont {Singh}},
  \bibinfo {author} {\bibfnamefont {A.~C.}\ \bibnamefont {Mishra}}, \ and\
  \bibinfo {author} {\bibfnamefont {S.~K.}\ \bibnamefont {Mohanty}},\
  }\bibfield  {title} {\enquote {\bibinfo {title} {Effects of both diffuse and
  collimated incident radiation on phototactic bioconvection},}\ }\href@noop {}
  {\bibfield  {journal} {\bibinfo  {journal} {Physics of Fluids}\ }\textbf
  {\bibinfo {volume} {28}},\ \bibinfo {pages} {124104} (\bibinfo {year}
  {2016})}\BibitemShut {NoStop}%
\bibitem [{\citenamefont {Panda}(2020)}]{8panda2020}%
  \BibitemOpen
  \bibfield  {author} {\bibinfo {author} {\bibfnamefont {M.~K.}\ \bibnamefont
  {Panda}},\ }\bibfield  {title} {\enquote {\bibinfo {title} {Effects of
  anisotropic scattering on the onset of phototactic bioconvection with diffuse
  and collimated irradiation},}\ }\href@noop {} {\bibfield  {journal} {\bibinfo
   {journal} {Physics of Fluids}\ }\textbf {\bibinfo {volume} {32}},\ \bibinfo
  {pages} {091903} (\bibinfo {year} {2020})}\BibitemShut {NoStop}%
\bibitem [{\citenamefont {Straughan}(1993)}]{9straughan1993}%
  \BibitemOpen
  \bibfield  {author} {\bibinfo {author} {\bibfnamefont {B.}~\bibnamefont
  {Straughan}},\ }\href@noop {} {\emph {\bibinfo {title} {Mathematical aspects
  of penetrative convection}}}\ (\bibinfo  {publisher} {CRC Press},\ \bibinfo
  {year} {1993})\BibitemShut {NoStop}%
\bibitem [{\citenamefont {Vincent}\ and\ \citenamefont
  {Hill}(1996)}]{12vincent1996}%
  \BibitemOpen
  \bibfield  {author} {\bibinfo {author} {\bibfnamefont {R.~V.}\ \bibnamefont
  {Vincent}}\ and\ \bibinfo {author} {\bibfnamefont {N.~A.}\ \bibnamefont
  {Hill}},\ }\bibfield  {title} {\enquote {\bibinfo {title} {Bioconvection in a
  suspension of phototactic algae},}\ }\href@noop {} {\bibfield  {journal}
  {\bibinfo  {journal} {Journal of Fluid Mechanics}\ }\textbf {\bibinfo
  {volume} {327}},\ \bibinfo {pages} {343--371} (\bibinfo {year}
  {1996})}\BibitemShut {NoStop}%
\bibitem [{\citenamefont {Ghorai}\ and\ \citenamefont
  {Hill}(2005)}]{10ghorai2005}%
  \BibitemOpen
  \bibfield  {author} {\bibinfo {author} {\bibfnamefont {S.}~\bibnamefont
  {Ghorai}}\ and\ \bibinfo {author} {\bibfnamefont {N.~A.}\ \bibnamefont
  {Hill}},\ }\bibfield  {title} {\enquote {\bibinfo {title} {Penetrative
  phototactic bioconvection},}\ }\href@noop {} {\bibfield  {journal} {\bibinfo
  {journal} {Physics of fluids}\ }\textbf {\bibinfo {volume} {17}},\ \bibinfo
  {pages} {074101} (\bibinfo {year} {2005})}\BibitemShut {NoStop}%
\bibitem [{\citenamefont {Ghorai}, \citenamefont {Panda},\ and\ \citenamefont
  {Hill}(2010)}]{7ghorai2010}%
  \BibitemOpen
  \bibfield  {author} {\bibinfo {author} {\bibfnamefont {S.}~\bibnamefont
  {Ghorai}}, \bibinfo {author} {\bibfnamefont {M.~K.}\ \bibnamefont {Panda}}, \
  and\ \bibinfo {author} {\bibfnamefont {N.~A.}\ \bibnamefont {Hill}},\
  }\bibfield  {title} {\enquote {\bibinfo {title} {Bioconvection in a
  suspension of isotropically scattering phototactic algae},}\ }\href@noop {}
  {\bibfield  {journal} {\bibinfo  {journal} {Physics of Fluids}\ }\textbf
  {\bibinfo {volume} {22}},\ \bibinfo {pages} {071901} (\bibinfo {year}
  {2010})}\BibitemShut {NoStop}%
\bibitem [{\citenamefont {Ghorai}\ and\ \citenamefont
  {Panda}(2013)}]{13ghorai2013}%
  \BibitemOpen
  \bibfield  {author} {\bibinfo {author} {\bibfnamefont {S.}~\bibnamefont
  {Ghorai}}\ and\ \bibinfo {author} {\bibfnamefont {M.~K.}\ \bibnamefont
  {Panda}},\ }\bibfield  {title} {\enquote {\bibinfo {title} {Bioconvection in
  an anisotropic scattering suspension of phototactic algae},}\ }\href@noop {}
  {\bibfield  {journal} {\bibinfo  {journal} {European Journal of
  Mechanics-B/Fluids}\ }\textbf {\bibinfo {volume} {41}},\ \bibinfo {pages}
  {81--93} (\bibinfo {year} {2013})}\BibitemShut {NoStop}%
\bibitem [{\citenamefont {Panda}\ and\ \citenamefont
  {Singh}(2016)}]{11panda2016}%
  \BibitemOpen
  \bibfield  {author} {\bibinfo {author} {\bibfnamefont {M.~K.}\ \bibnamefont
  {Panda}}\ and\ \bibinfo {author} {\bibfnamefont {R.}~\bibnamefont {Singh}},\
  }\bibfield  {title} {\enquote {\bibinfo {title} {Penetrative phototactic
  bioconvection in a two-dimensional non-scattering suspension},}\ }\href@noop
  {} {\bibfield  {journal} {\bibinfo  {journal} {Physics of Fluids}\ }\textbf
  {\bibinfo {volume} {28}},\ \bibinfo {pages} {054105} (\bibinfo {year}
  {2016})}\BibitemShut {NoStop}%
\bibitem [{\citenamefont {Panda}, \citenamefont {Sharma},\ and\ \citenamefont
  {Kumar}(2022)}]{16panda2022}%
  \BibitemOpen
  \bibfield  {author} {\bibinfo {author} {\bibfnamefont {M.~K.}\ \bibnamefont
  {Panda}}, \bibinfo {author} {\bibfnamefont {P.}~\bibnamefont {Sharma}}, \
  and\ \bibinfo {author} {\bibfnamefont {S.}~\bibnamefont {Kumar}},\ }\bibfield
   {title} {\enquote {\bibinfo {title} {Effects of oblique irradiation on the
  onset of phototactic bioconvection},}\ }\href@noop {} {\bibfield  {journal}
  {\bibinfo  {journal} {Physics of Fluids}\ }\textbf {\bibinfo {volume} {34}},\
  \bibinfo {pages} {024108} (\bibinfo {year} {2022})}\BibitemShut {NoStop}%
\bibitem [{\citenamefont {Kumar}(2022)}]{17kumar2022}%
  \BibitemOpen
  \bibfield  {author} {\bibinfo {author} {\bibfnamefont {S.}~\bibnamefont
  {Kumar}},\ }\bibfield  {title} {\enquote {\bibinfo {title} {Phototactic
  isotropic scattering bioconvection with oblique irradiation},}\ }\href@noop
  {} {\bibfield  {journal} {\bibinfo  {journal} {Physics of Fluids}\ }\textbf
  {\bibinfo {volume} {34}},\ \bibinfo {pages} {114125} (\bibinfo {year}
  {2022})}\BibitemShut {NoStop}%
\bibitem [{\citenamefont {Kumar}(2023)}]{39kumar2023}%
  \BibitemOpen
  \bibfield  {author} {\bibinfo {author} {\bibfnamefont {S.}~\bibnamefont
  {Kumar}},\ }\bibfield  {title} {\enquote {\bibinfo {title} {Isotropic
  scattering with a rigid upper surface at the onset of phototactic
  bioconvection},}\ }\href@noop {} {\bibfield  {journal} {\bibinfo  {journal}
  {Physics of Fluids}\ }\textbf {\bibinfo {volume} {35}},\ \bibinfo {pages}
  {024106} (\bibinfo {year} {2023})}\BibitemShut {NoStop}%
\bibitem [{\citenamefont {Cash}\ and\ \citenamefont
  {Moore}(1980)}]{19cash1980}%
  \BibitemOpen
  \bibfield  {author} {\bibinfo {author} {\bibfnamefont {J.~R.}\ \bibnamefont
  {Cash}}\ and\ \bibinfo {author} {\bibfnamefont {D.~R.}\ \bibnamefont
  {Moore}},\ }\bibfield  {title} {\enquote {\bibinfo {title} {A high order
  method for the numerical solution of two-point boundary value problems},}\
  }\href@noop {} {\bibfield  {journal} {\bibinfo  {journal} {BIT Numerical
  Mathematics}\ }\textbf {\bibinfo {volume} {20}},\ \bibinfo {pages} {44--52}
  (\bibinfo {year} {1980})}\BibitemShut {NoStop}%
\bibitem [{\citenamefont {Kitsunezaki}, \citenamefont {Komori},\ and\
  \citenamefont {Harumoto}(2007)}]{2kitsunezaki2007}%
  \BibitemOpen
  \bibfield  {author} {\bibinfo {author} {\bibfnamefont {S.}~\bibnamefont
  {Kitsunezaki}}, \bibinfo {author} {\bibfnamefont {R.}~\bibnamefont {Komori}},
  \ and\ \bibinfo {author} {\bibfnamefont {T.}~\bibnamefont {Harumoto}},\
  }\bibfield  {title} {\enquote {\bibinfo {title} {Bioconvection and front
  formation of paramecium tetraurelia},}\ }\href@noop {} {\bibfield  {journal}
  {\bibinfo  {journal} {Physical Review E}\ }\textbf {\bibinfo {volume} {76}},\
  \bibinfo {pages} {046301} (\bibinfo {year} {2007})}\BibitemShut {NoStop}%
\bibitem [{\citenamefont {Kessler}(1989)}]{5kessler1989}%
  \BibitemOpen
  \bibfield  {author} {\bibinfo {author} {\bibfnamefont {J.}~\bibnamefont
  {Kessler}},\ }\bibfield  {title} {\enquote {\bibinfo {title} {Path and
  pattern-the mutual dynamics of swimming cells and their environment},}\
  }\href@noop {} {\bibfield  {journal} {\bibinfo  {journal} {Comments Theor.
  Biol.}\ }\textbf {\bibinfo {volume} {1}},\ \bibinfo {pages} {85--108}
  (\bibinfo {year} {1989})}\BibitemShut {NoStop}%
\bibitem [{\citenamefont {Panda}\ and\ \citenamefont
  {Ghorai}(2013)}]{14panda2013}%
  \BibitemOpen
  \bibfield  {author} {\bibinfo {author} {\bibfnamefont {M.~K.}\ \bibnamefont
  {Panda}}\ and\ \bibinfo {author} {\bibfnamefont {S.}~\bibnamefont {Ghorai}},\
  }\bibfield  {title} {\enquote {\bibinfo {title} {Penetrative phototactic
  bioconvection in an isotropic scattering suspension},}\ }\href@noop {}
  {\bibfield  {journal} {\bibinfo  {journal} {Physics of Fluids}\ }\textbf
  {\bibinfo {volume} {25}},\ \bibinfo {pages} {071902} (\bibinfo {year}
  {2013})}\BibitemShut {NoStop}%
\bibitem [{\citenamefont {Hill}\ and\ \citenamefont
  {H{\"a}der}(1997)}]{18hill1997}%
  \BibitemOpen
  \bibfield  {author} {\bibinfo {author} {\bibfnamefont {N.~A.}\ \bibnamefont
  {Hill}}\ and\ \bibinfo {author} {\bibfnamefont {D.-P.}\ \bibnamefont
  {H{\"a}der}},\ }\bibfield  {title} {\enquote {\bibinfo {title} {A biased
  random walk model for the trajectories of swimming micro-organisms},}\
  }\href@noop {} {\bibfield  {journal} {\bibinfo  {journal} {Journal of
  theoretical biology}\ }\textbf {\bibinfo {volume} {186}},\ \bibinfo {pages}
  {503--526} (\bibinfo {year} {1997})}\BibitemShut {NoStop}%
\bibitem [{\citenamefont {Kage}\ \emph {et~al.}(2013)\citenamefont {Kage},
  \citenamefont {Hosoya}, \citenamefont {Baba},\ and\ \citenamefont
  {Mogami}}]{27kage2013}%
  \BibitemOpen
  \bibfield  {author} {\bibinfo {author} {\bibfnamefont {A.}~\bibnamefont
  {Kage}}, \bibinfo {author} {\bibfnamefont {C.}~\bibnamefont {Hosoya}},
  \bibinfo {author} {\bibfnamefont {S.~A.}\ \bibnamefont {Baba}}, \ and\
  \bibinfo {author} {\bibfnamefont {Y.}~\bibnamefont {Mogami}},\ }\bibfield
  {title} {\enquote {\bibinfo {title} {Drastic reorganization of the
  bioconvection pattern of chlamydomonas: quantitative analysis of the pattern
  transition response},}\ }\href@noop {} {\bibfield  {journal} {\bibinfo
  {journal} {Journal of Experimental Biology}\ }\textbf {\bibinfo {volume}
  {216}},\ \bibinfo {pages} {4557--4566} (\bibinfo {year} {2013})}\BibitemShut
  {NoStop}%
\bibitem [{\citenamefont {Mendelson}\ and\ \citenamefont
  {Lega}(1998)}]{28mendelson1998}%
  \BibitemOpen
  \bibfield  {author} {\bibinfo {author} {\bibfnamefont {N.~H.}\ \bibnamefont
  {Mendelson}}\ and\ \bibinfo {author} {\bibfnamefont {J.}~\bibnamefont
  {Lega}},\ }\bibfield  {title} {\enquote {\bibinfo {title} {A complex pattern
  of traveling stripes is produced by swimming cells of bacillus subtilis},}\
  }\href@noop {} {\bibfield  {journal} {\bibinfo  {journal} {Journal of
  bacteriology}\ }\textbf {\bibinfo {volume} {180}},\ \bibinfo {pages}
  {3285--3294} (\bibinfo {year} {1998})}\BibitemShut {NoStop}%
\bibitem [{\citenamefont {Gittleson}\ and\ \citenamefont
  {Jahn}(1968)}]{29gittleson1968}%
  \BibitemOpen
  \bibfield  {author} {\bibinfo {author} {\bibfnamefont {S.~M.}\ \bibnamefont
  {Gittleson}}\ and\ \bibinfo {author} {\bibfnamefont {T.~L.}\ \bibnamefont
  {Jahn}},\ }\bibfield  {title} {\enquote {\bibinfo {title} {Pattern swimming
  by polytomella agilis},}\ }\href@noop {} {\bibfield  {journal} {\bibinfo
  {journal} {The American Naturalist}\ }\textbf {\bibinfo {volume} {102}},\
  \bibinfo {pages} {413--425} (\bibinfo {year} {1968})}\BibitemShut {NoStop}%
\bibitem [{\citenamefont {Khan}\ \emph {et~al.}(2017)\citenamefont {Khan},
  \citenamefont {Gul}, \citenamefont {Khan}, \citenamefont {Bonyah},\ and\
  \citenamefont {Islam}}]{30khan2017}%
  \BibitemOpen
  \bibfield  {author} {\bibinfo {author} {\bibfnamefont {N.~S.}\ \bibnamefont
  {Khan}}, \bibinfo {author} {\bibfnamefont {T.}~\bibnamefont {Gul}}, \bibinfo
  {author} {\bibfnamefont {M.~A.}\ \bibnamefont {Khan}}, \bibinfo {author}
  {\bibfnamefont {E.}~\bibnamefont {Bonyah}}, \ and\ \bibinfo {author}
  {\bibfnamefont {S.}~\bibnamefont {Islam}},\ }\bibfield  {title} {\enquote
  {\bibinfo {title} {Mixed convection in gravity-driven thin film non-newtonian
  nanofluids flow with gyrotactic microorganisms},}\ }\href@noop {} {\bibfield
  {journal} {\bibinfo  {journal} {Results in physics}\ }\textbf {\bibinfo
  {volume} {7}},\ \bibinfo {pages} {4033--4049} (\bibinfo {year}
  {2017})}\BibitemShut {NoStop}%
\bibitem [{\citenamefont {Hayat}, \citenamefont {Alsaedi}\ \emph
  {et~al.}(2021)\citenamefont {Hayat}, \citenamefont {Alsaedi} \emph
  {et~al.}}]{31hayat2021}%
  \BibitemOpen
  \bibfield  {author} {\bibinfo {author} {\bibfnamefont {T.}~\bibnamefont
  {Hayat}}, \bibinfo {author} {\bibfnamefont {A.}~\bibnamefont {Alsaedi}},
  \emph {et~al.},\ }\bibfield  {title} {\enquote {\bibinfo {title} {Development
  of bioconvection flow of nanomaterial with melting effects},}\ }\href@noop {}
  {\bibfield  {journal} {\bibinfo  {journal} {Chaos, Solitons \& Fractals}\
  }\textbf {\bibinfo {volume} {148}},\ \bibinfo {pages} {111015} (\bibinfo
  {year} {2021})}\BibitemShut {NoStop}%
\bibitem [{\citenamefont {Incropera}, \citenamefont {Wagner},\ and\
  \citenamefont {Houf}(1981)}]{32incropera1981}%
  \BibitemOpen
  \bibfield  {author} {\bibinfo {author} {\bibfnamefont {F.~P.}\ \bibnamefont
  {Incropera}}, \bibinfo {author} {\bibfnamefont {T.~R.}\ \bibnamefont
  {Wagner}}, \ and\ \bibinfo {author} {\bibfnamefont {W.~G.}\ \bibnamefont
  {Houf}},\ }\bibfield  {title} {\enquote {\bibinfo {title} {A comparison of
  predictions and measurements of the radiation field in a shallow water
  layer},}\ }\href@noop {} {\bibfield  {journal} {\bibinfo  {journal} {Water
  Resources Research}\ }\textbf {\bibinfo {volume} {17}},\ \bibinfo {pages}
  {142--148} (\bibinfo {year} {1981})}\BibitemShut {NoStop}%
\bibitem [{\citenamefont {Daniel}, \citenamefont {Laurendeau},\ and\
  \citenamefont {Incropera}(1979)}]{33daniel1979}%
  \BibitemOpen
  \bibfield  {author} {\bibinfo {author} {\bibfnamefont {K.~J.}\ \bibnamefont
  {Daniel}}, \bibinfo {author} {\bibfnamefont {N.~M.}\ \bibnamefont
  {Laurendeau}}, \ and\ \bibinfo {author} {\bibfnamefont {F.~P.}\ \bibnamefont
  {Incropera}},\ }\bibfield  {title} {\enquote {\bibinfo {title} {Prediction of
  radiation absorption and scattering in turbid water bodies},}\ }\href@noop {}
  {\  (\bibinfo {year} {1979})}\BibitemShut {NoStop}%
\bibitem [{\citenamefont {Hill}, \citenamefont {Pedley},\ and\ \citenamefont
  {Kessler}(1989)}]{34hill1989}%
  \BibitemOpen
  \bibfield  {author} {\bibinfo {author} {\bibfnamefont {N.~A.}\ \bibnamefont
  {Hill}}, \bibinfo {author} {\bibfnamefont {T.~J.}\ \bibnamefont {Pedley}}, \
  and\ \bibinfo {author} {\bibfnamefont {J.~O.}\ \bibnamefont {Kessler}},\
  }\bibfield  {title} {\enquote {\bibinfo {title} {Growth of bioconvection
  patterns in a suspension of gyrotactic micro-organisms in a layer of finite
  depth},}\ }\href@noop {} {\bibfield  {journal} {\bibinfo  {journal} {Journal
  of Fluid Mechanics}\ }\textbf {\bibinfo {volume} {208}},\ \bibinfo {pages}
  {509--543} (\bibinfo {year} {1989})}\BibitemShut {NoStop}%
\bibitem [{\citenamefont {Modest}\ and\ \citenamefont
  {Mazumder}(2021)}]{35modest2021}%
  \BibitemOpen
  \bibfield  {author} {\bibinfo {author} {\bibfnamefont {M.~F.}\ \bibnamefont
  {Modest}}\ and\ \bibinfo {author} {\bibfnamefont {S.}~\bibnamefont
  {Mazumder}},\ }\href@noop {} {\emph {\bibinfo {title} {Radiative heat
  transfer}}}\ (\bibinfo  {publisher} {Academic press},\ \bibinfo {year}
  {2021})\BibitemShut {NoStop}%
\bibitem [{\citenamefont {Chandrasekhar}(1960)}]{36chandrasekhar1960}%
  \BibitemOpen
  \bibfield  {author} {\bibinfo {author} {\bibfnamefont {S.}~\bibnamefont
  {Chandrasekhar}},\ }\bibfield  {title} {\enquote {\bibinfo {title} {Radiative
  transfer dover publications inc},}\ }\href@noop {} {\bibfield  {journal}
  {\bibinfo  {journal} {New York}\ } (\bibinfo {year} {1960})}\BibitemShut
  {NoStop}%
\bibitem [{\citenamefont {Ghorai}\ and\ \citenamefont
  {Singh}(2009)}]{37ghorai2009}%
  \BibitemOpen
  \bibfield  {author} {\bibinfo {author} {\bibfnamefont {S.}~\bibnamefont
  {Ghorai}}\ and\ \bibinfo {author} {\bibfnamefont {R.}~\bibnamefont {Singh}},\
  }\bibfield  {title} {\enquote {\bibinfo {title} {Linear stability analysis of
  gyrotactic plumes},}\ }\href@noop {} {\bibfield  {journal} {\bibinfo
  {journal} {Physics of Fluids}\ }\textbf {\bibinfo {volume} {21}},\ \bibinfo
  {pages} {081901} (\bibinfo {year} {2009})}\BibitemShut {NoStop}%
\bibitem [{\citenamefont {Press}(1992)}]{38press1992}%
  \BibitemOpen
  \bibfield  {author} {\bibinfo {author} {\bibfnamefont {W.~H.}\ \bibnamefont
  {Press}},\ }\bibfield  {title} {\enquote {\bibinfo {title} {Numerical recipes
  in fortran.}}\ }\href@noop {} {\bibfield  {journal} {\bibinfo  {journal} {The
  Art of Scientific Computing.}\ } (\bibinfo {year} {1992})}\BibitemShut
  {NoStop}%
\end{thebibliography}%
	
\end{document}